\documentclass{amsproc}

\usepackage{graphicx}
\usepackage{amsfonts}
\usepackage{amscd}
\usepackage{amssymb}
\usepackage[all]{xypic}
\usepackage{tikz}

\newtheorem{theorem}{Theorem}[section]
\newtheorem{proposition}[theorem]{Proposition}

\newtheorem{corollary}[theorem]{Corollary}

\theoremstyle{definition}
\newtheorem{definition}[theorem]{Definition}
\newtheorem{Notation}[theorem]{Notation}
\newtheorem{example}[theorem]{Example}

\theoremstyle{remark}
\newtheorem{remark}[theorem]{Remark}
\newtheorem{question}[theorem]{Question}

\numberwithin{equation}{section}

\DeclareMathOperator{\Gal}{Gal}

\DeclareMathOperator{\sep}{sep}
\DeclareMathOperator{\val}{val}
\DeclareMathOperator{\tame}{tame}
\DeclareMathOperator{\sconn}{sc}
\DeclareMathOperator{\mult}{m}
\DeclareMathOperator{\Mult}{mult}
\DeclareMathOperator{\add}{a}
\DeclareMathOperator{\Zar}{Zar}

\DeclareMathOperator{\Hom}{Hom}

\DeclareMathOperator{\Id}{Id}
\DeclareMathOperator{\Ker}{Ker}
\DeclareMathOperator{\Tor}{Tor}
\DeclareMathOperator{\Bor}{Bor}
\DeclareMathOperator{\BorTor}{BorTor}
\DeclareMathOperator{\Spr}{Spr}
\DeclareMathOperator{\Norm}{Norm}
\DeclareMathOperator{\Int}{Int}
\DeclareMathOperator{\Pic}{Pic}
\DeclareMathOperator{\pr}{pr}
\DeclareMathOperator{\Kott}{Kott}
\DeclareMathOperator{\Ext}{Ext}
\DeclareMathOperator{\CExt}{CExt}
\DeclareMathOperator{\Comm}{Comm}
\DeclareMathOperator{\Sym}{Sym}
\DeclareMathOperator{\Aff}{Aff}
\DeclareMathOperator{\supremum}{sup}
\DeclareMathOperator{\characteristic}{char}
\DeclareMathOperator{\Leg}{Leg}

\newcommand{\Matrix}[4]{ \left( \begin{array}{cc}  #1 & #2 \\  #3 & #4 \\ \end{array} \right) }

\newcommand{\TMatrix}[9]{ \left( \begin{array}{ccc} #1 & #2 & #3 \\  #4 & #5 & #6 \\  #7 & #8 & #9  \\ \end{array} \right) }

\newcommand{\Type}[1]{\mathbf{\mathsf{#1}}}

\newcommand{\ideal}[1]{ {\mathfrak{#1}} }
\newcommand{\Lie}[1]{ {\mathfrak{#1}} }
\newcommand{\Cat}[1]{ {\mathfrak{#1}} }

\newcommand{\alg}[1]{{\boldsymbol{#1}} }

\newcommand{\ZZ}{\mathbb Z}
\newcommand{\FF}{\mathbb F}
\newcommand{\LL}{\mathbb L}
\newcommand{\QQ}{\mathbb Q}
\newcommand{\RR}{\mathbb R}

\newcommand{\HH}{\mathbb H}
\newcommand{\adeles}{\mathbb A}

\newcommand{\Build}{{\mathcal B}}
\newcommand{\Apart}{{\mathcal A}}
\newcommand{\OO}{\mathcal{O}}

\newcommand{\isom}{\cong}

\begin{document}

\title{Managing metaplectiphobia:  Covering p-adic groups.}

\author{Martin H. Weissman}
\address{Department of Mathematics, University of California, Santa Cruz, California 95064}
\email{weissman@ucsc.edu}

\date{July 31, 2010.}



\begin{abstract}
Brylinski and Deligne have provided a framework to study central extensions of reductive groups by $\alg{K}_2$ over a field $F$.  Such central extensions can be used to construct central extensions of $p$-adic groups by finite cyclic groups, including the metaplectic groups.  

Particularly interesting is the observation of Brylinski and Deligne that a central extension of a reductive group by $\alg{K}_2$, over a $p$-adic field, yields a family of central extensions of reductive groups by the multiplicative group over the residue field, indexed by the points of the building.  These algebraic groups over the residue field determine the structure of central extensions of $p$-adic groups, when the extension is restricted to a parahoric subgroup.  

This article surveys and builds upon the work of Brylinski and Deligne, culminating in a precise description of some central extensions using the Bruhat-Tits building.
\end{abstract}
\maketitle


\section*{Preliminaries}
\subsection{Introduction}
Let $F$ be a $p$-adic field, with ring of integers $\OO$ and residue field $\FF$.  Let $\alg{G}$ be a connected reductive group over $F$, and let $G = \alg{G}(F)$.  We are interested in a class of ``tame'' central extensions of $G$ by $\FF^\times$:
$$1 \rightarrow \FF^\times \rightarrow \tilde G \rightarrow G \rightarrow 1.$$
Some authors study all such central extensions in the category of locally compact topological groups; while this is certainly possible, we are compelled to work with a different (and effectively narrower) category of central extensions arising from a construction of Brylinski and Deligne \cite{B-D}.  They begin with a central extension in the category of sheaves of groups on the big Zariski site (a category which includes the category of algebraic groups as a full subcategory) over $F$:
$$1 \rightarrow \alg{K}_2 \rightarrow \alg{G}' \rightarrow \alg{G} \rightarrow 1.$$
From such a central extension, one may take $F$-points to get an exact sequence of groups $\alg{K}_2(F) \rightarrow \alg{G}'(F) \rightarrow G$, and push forward using the tame symbol in K-theory:
$$\tame: \alg{K}_2(F) \rightarrow \FF^\times.$$
This yields an extension of locally compact topological groups:
$$1 \rightarrow \FF^\times \rightarrow \tilde G \rightarrow G \rightarrow 1.$$

There are many reasons for considering only extensions arising through this construction of Brylinski and Deligne, rather than a more general class of topological central extensions.  We list some reasons below:
\begin{enumerate}
\item
When considering Brylinski-Deligne central extensions globally, the central extension of the adelic group splits  canonically over the rational points of the group, leading to a reasonable definition of automorphic forms and representations.  Here it must be mentioned that Prasad-Raghunathan \cite{PRag} and Prasad-Rapinchuk \cite{PRap} have determined the ``metaplectic kernel,'' which in turn describes all central extensions of $\alg{G}(\adeles_F)$ by finite abelian groups which split canonically over $\alg{G}(F)$, when $F$ is a global field and $\alg{G}$ is absolutely simple and simply connected over $F$.  Thus central extensions constructed by Brylinski-Deligne (over a global field and its adeles) fit into a class of metaplectic groups studied by Prasad, Raghunathan, Rapinchuk, and others.
\item
Brylinski and Deligne have classified their central extensions by essentially combinatorial data related to the root datum.  If one hopes for a Langlands-style conjecture for central extensions, then one must have such a combinatorial classification to speculate about an analogue of the ``Langlands dual group''.
\item
The work of Brylinski and Deligne does not rely on choosing specific cocycles; their work describes a {\em category} of central extensions, not only the isomorphism classes thereof.  This is crucial, since any putative parameterization of representations of $\tilde G$ would depend on the choice of cocycle.  This arises in practice, where parameterizations of representations of metaplectic groups and descriptions of Hecke algebras of metaplectic groups depend on initial choices (usually a choice of signs (choosing $i$ or $-i$) or cocycles).
\item
Brylinski and Deligne demonstrate a remarkable connection between central extensions of reductive groups by $\alg{K}_2$, over $F$, and central extensions of related (by Bruhat-Tits theory) reductive groups by $\alg{G}_{\mult}$, over $\FF$.  Thus ultimately (from the standpoint of K-types) the representation theory of $\tilde G$ boils down to representation theory of finite groups of Lie type.
\end{enumerate}

It is this last point which is the focus of this article.  When $\tilde G$ is a central extension of $G$ by $\FF^\times$, obtained from the construction of Brylinski and Deligne, and $x$ is a point in the Bruhat-Tits building of $G$, one may restrict the central extension to the parahoric subgroup $G_x$ to obtain:
\begin{equation}
\label{ParahoricExtension}
1 \rightarrow \FF^\times \rightarrow \tilde G_x \rightarrow G_x \rightarrow 1.
\end{equation}

By Bruhat-Tits theory, the quotient $\bar M_x = G_x / G_x^+$ of the parahoric $G_x$ by a pro-$p$ subgroup $G_x^+$ coincides with the $\FF$ points of a connected reductive group $\alg{\bar M}_x$ over $\FF$.  The central extension $\tilde G_x$ splits canonically over $G_x^+$, leading to a central extension of finite groups:
$$1 \rightarrow \FF^\times \rightarrow \bar M_x' \rightarrow \bar M_x \rightarrow 1.$$
In Construction 12.11 of \cite{B-D}, it is shown that this central extension of finite groups arises from a central extension of algebraic groups over $\FF$:
\begin{equation}
\label{ResidualExtension}
1 \rightarrow \alg{\bar G}_{\mult} \rightarrow \alg{\bar M}_x' \rightarrow \alg{\bar M}_x \rightarrow 1.
\end{equation}
The central extension (\ref{ParahoricExtension}) is uniquely determined by the central extension (\ref{ResidualExtension}) of the reductive group $\alg{\bar M}_x$ by $\alg{\bar G}_{\mult}$ over the residue field.

Deligne and Brylinski ask a natural question, listed as Question 12.13(i) of \cite{B-D}, and directly quoted below:
\begin{quote}
Suppose that $G$ is reductive, and that $E$ is given as in 7.2, for $T$ a maximally split maximal torus of $G$.  Suppose that $G_V$ is given as in Bruhat-Tits (1984) 4.6.  It would be interesting to compute the central extension $G_s^\sim$ in that case, especially for $G_V(V)$ a maximal bounded subgroup of $G(K)$, given by a vertex of the building of $G$.
\end{quote}
Rephrased in our notation, Deligne and Brylinski ask:
\begin{question}
\label{Q2}
What is the central extension $\alg{\bar M}_x'$ of $\alg{\bar M}_x$ by $\alg{\bar G}_{\mult}$?
\end{question}

For better or worse, answering this question requires an answer to a general question about reductive groups over fields:
\begin{question}
\label{Q3}
Given a connected reductive group $\alg{G}$ over a field $F$, describe the category of central extensions of $\alg{G}$ by $\alg{G}_{\mult}$ over $F$.
\end{question}
If we wish to describe the central extension $\alg{\bar M}_x'$ of $\alg{\bar M}_x$ by $\alg{\bar G}_{\mult}$ (Question \ref{Q2}), we require a general method of describing such central extensions (Question \ref{Q3}).  One may describe $\alg{\bar M}_x'$ up to isomorphism, by describing the root datum of $\alg{\bar M}_x'$ along with maps of cocharacter lattices $\ZZ \rightarrow Y' \rightarrow Y$ corresponding to the central extension $\alg{\bar G}_{\mult} \rightarrow \alg{\bar M}_x' \rightarrow \alg{M}_x$.  But, for reasons of descent and non-ambiguous parameterization of representations, this is insufficient.  One must go further and describe $\alg{\bar M}_x'$ up to {\em unique} isomorphism.  This description is new, and is given by our Theorem \ref{CEByGm}.  It should be said that our Theorem \ref{CEByGm} is analogous to (but does not follow from) the Main Theorem of \cite{B-D} -- essentially we consider central extensions of reductive groups by $\alg{K}_1$ while Brylinski and Deligne consider central extensions of reductive groups by $\alg{K}_2$.   

After we answer Question \ref{Q3} with Theorem \ref{CEByGm}, we are able to answer Question \ref{Q2} to a large extent.  Without providing a general answer to Question \ref{Q2}, we provide the necessary tools, and illustrate this with some examples.  Our examples include the simplest case $\alg{G} = \alg{SL}_2$ to illustrate some basic principles, $\alg{SU}_3$ to illustrate the non-split case, and $\alg{G} = \alg{G}_2$ to demonstrate how our methods generalize.  Many of our calculations, in the quasisplit case, have been carried out by Deodhar \cite{Deo}, in a somewhat different framework, and with different goals in mind.

By answering Question \ref{Q2}, we are able to describe covers of parahoric subgroups in tame central extensions of $p$-adic groups.  This complements earlier work \cite{We2} (joint with T. Howard) on depth zero representations of these central extensions.  Indeed an answer to Question \ref{Q2} seems crucial, if one wishes to find an appropriate generalization of the (local) Langlands conjectures to nonlinear covering groups.

Beyond answering Questions \ref{Q2} and \ref{Q3}, we hope that this article serves as a guide for others who wish to use Brylinski and Deligne's framework when studying metaplectic groups and more general nonlinear covers of $p$-adic groups.  At the very least, we hope to demonstrate the strength and elegance of \cite{B-D}, by surveying and expanding upon their results. 

\subsection{Notation}
$F$ will always denote a field, with a separable closure $F^{\sep}$.  We use a boldface font, like $\alg{J}$ for an algebraic variety over $F$, or more generally for any functor from the category of finitely-generated $F$-algebras to the category of sets.  When $A$ is a finitely-generated $F$-algebra, we write $\alg{J}(A)$ for the $A$-points of $\alg{J}$; more generally, for any $F$-algebra $A$, we write $\alg{J}(A)$ for the direct limit of the $A_i$-points of $\alg{J}$, as $A_i$ ranges over the directed set of finitely-generated sub-$F$-algebras of $A$.  We use an ordinary font for the $F$-points:  $J = \alg{J}(F)$.  Similarly, we use a boldface font, like $\alg{j}: \alg{J} \rightarrow \alg{K}$ for a morphism of algebraic varieties over $F$, or more generally for a natural transformation of set-valued functors on the category of finitely-generated $F$-algebras.  We use an ordinary font for the resulting function on $F$-points, as in $j: J \rightarrow K$.  When defining a morphism $\alg{j}$, we often just describe the function $j$ on $F$-points, leaving it to the reader to infer its algebraic origin.

When $\alg{p}: \alg{G}' \rightarrow \alg{G}$ is a surjective homomorphism of groups over $F$, a {\em section} of $\alg{p}$ means an algebraic map $\alg{j}: \alg{G} \rightarrow \alg{G}'$ satisfying $\alg{p} \circ \alg{j} = \alg{\Id}$.  A {\em splitting} of $\alg{p}$ is a section which is also a homomorphism.  If $\alg{\Ker}(\alg{p})$ is central in $\alg{G}'$, and $\alg{\chi}: \alg{G} \rightarrow \alg{\Ker}(\alg{p})$ is a homomorphism, then we may {\em twist} a section or splitting $\alg{j}$ by $\alg{\chi}$: $\alg{j} \cdot \alg{\chi}$ is also a section or splitting, accordingly.  We use similar terminology, for surjective homomorphisms of abstract groups (using abstract maps and homomorphisms), and group-valued functors (using natural transformations of set-valued functors, and natural transformations of group-valued functors).

Eventually, we will assume that $F$ is a field with nontrivial discrete valuation $\val: F^\times \rightarrow \RR$.  In this circumstance we let $\OO$ be the valuation ring of $F$, and $\ideal{p}$ the maximal ideal of $\OO$.  The residue field $\FF = \OO / \ideal{p}$ will always be assumed perfect.  We use an overline when working over $\FF$; thus $\alg{\bar J}$ might denote an algebraic variety over $\FF$, and $\bar J$ its $\FF$-points.  We use an underline when working over $\OO$; thus $\alg{\underline J}$ might denote a scheme over $\OO$; in this situation, we would write $\alg{J}$ for its generic fibre -- a scheme over $F$ -- and $\alg{\bar J}$ for its special fibre -- a scheme over $\FF$.  We follow this convention also for morphisms:  $\alg{\underline j}$ might denote a morphism of schemes over $\OO$, and $\alg{\bar j}$ a morphism of schemes over $\FF$.  

The letter $\alg{G}$ will always denote an affine algebraic group over a field $F$.  We always write $\alg{G}_{\mult}$ for the multiplicative group over $F$, and $\alg{G}_{\add}$ for the additive group over $F$.  When $\alg{S}$ is a torus over a field $F$, we define $X(\alg{S}) = \Hom(\alg{S}, \alg{G}_{\mult})$ and $Y(\alg{S}) = \Hom(\alg{G}_{\mult}, \alg{S})$.  These are viewed as \'etale sheaves over $F$, or simply as abelian groups with an action of $\Gal(F^{\sep} / F)$.  

\subsection{K-groups}

For $n \geq 0$, we write $\alg{K}_n$ for the K-theory functor, from the category of finitely-generated $F$-algebras to the category of abelian groups.  We will thankfully only require reference to $\alg{K}_0$, $\alg{K}_1$, and $\alg{K}_2$ in this article.  We also will only require calculations of these groups for very simple classes of $F$-algebras.  Later we will view $\alg{K}_n$ as sheaves on the big Zariski site of schemes of finite type over $F$.

Whenever $A$ is a ring, $\alg{K}_0(A)$ is the Grothendieck group of finitely-generated projective $A$-modules.  In particular, whenever $L$ is a field, $\alg{K}_0(L)$ is naturally isomorphic to $\ZZ$, sending a finite-dimensional $L$-vector space to its dimension.  

Whenever $A$ is a Euclidean domain, $\alg{K}_1(A) = \alg{G}_{\mult}(A) = A^\times$.  In particular, when $L$ is a field, $\alg{K}_1(L) = L^\times$.    

It is somewhat difficult to define $\alg{K}_2(A)$ when $A$ is not a field.  However, for fields $L$ we have the following description: 
$$\alg{K}_2(L) = \frac{ L^\times \otimes_\ZZ L^\times }{ \langle x \otimes (1-x) \rangle_{1 \neq x \in L^\times} }.$$
When $x,y \in L^\times$, we write $\{ x, y \}$ for the image of $x \otimes y$ in $\alg{K}_2(L)$.  This {\em Steinberg symbol} satsfies the following relations:
\begin{description}
\item[Bilinearity]
$\{ x x', y \} = \{ x,y \} \{ x', y \}$ and $\{ x,yy' \} = \{ x,y \} \{ x, y' \}$ for all $x,x',y,y' \in L^\times$.
\item[Steinberg relation]
$\{x, 1-x \} = 1$ for all $1 \neq x \in L^\times$.
\item[Skew-symmetry]
$\{ x,y \} \{y, x \} = 1$ for all $x,y \in L^\times$.
\end{description}
In fact, skew-symmetry follows from the previous two properties.  Steinberg symbols are often not alternating, but they do satisfy the properties:  $\{ x, -x \} = 1$ and $\{x, x \} = \{ x, -1 \}$, for all $x \in L^\times$.  The group $\alg{K}_2(L)$ can be viewed as the abelian group generated by all formal symbols $\{x,y \}$ for $x,y \in L^\times$, modulo the relations above.  

\subsection{Acknowledgments}

We thank Brian Conrad and Mikhail Borovoi for providing some helpful references about algebraic groups and algebraic geometry.  In addition, we thank the anonymous referee for providing corrections and excellent suggestions to strengthen the exposition.  We thank the organizers, Loren Spice, Robert Doran, and Paul J. Sally, Jr., for inviting this paper.

\section{Central extensions by $\alg{G}_{\mult}$}
In this section, we let $F$ be a perfect field.  Let $\alg{G}$ be a connected reductive group over $F$.
\begin{definition}
A {\em central extension} of $\alg{G}$ by $\alg{G}_{\mult}$ is a triple $(\alg{G}', \alg{p}, \alg{\iota})$ where $\alg{G}'$ is an algebraic group over $F$, and $\alg{p},\alg{\iota}$ are morphisms of groups over $F$ fitting into a short exact sequence:
$$\xymatrix{
1 \ar[r] & \alg{G}_{\mult} \ar[r]^{\alg{\iota}} & \alg{G}' \ar[r]^{\alg{p}} & \alg{G} \ar[r] & 1,}$$
such that $\alg{\iota}$ is a closed embedding of $\alg{G}_{\mult}$ into the center of $\alg{G}'$, and $\alg{p}$ identifies $\alg{G}$ with the quotient group $\alg{G}' / \alg{\iota}(\alg{G}_{\mult})$.  
\end{definition}
Given a central extension $(\alg{G}',\alg{p},\alg{\iota})$ of $\alg{G}$ by $\alg{G}_{\mult}$, and any field $L$ containing $F$, Hilbert's Theorem 90 gives a short exact sequence of groups:
$$1 \rightarrow L^\times \rightarrow \alg{G}'(L) \rightarrow \alg{G}(L) \rightarrow 1.$$
When we write $g'$ for an element of $\alg{G}'(L)$, we always write $g$ for the projection, $g = p(g')$, in $\alg{G}(L)$.

\subsection{The category of central extensions}

\begin{definition}
Let $(\alg{G}_1', \alg{p}_1, \alg{\iota}_1)$ and $(\alg{G}_2', \alg{p}_2, \alg{\iota}_2)$ be two central extensions of $\alg{G}$ by $\alg{G}_{\mult}$.  A {\em morphism} from $(\alg{G}_1', \alg{p}_1, \alg{\iota}_1)$ to $(\alg{G}_2', \alg{p}_2, \alg{\iota}_2)$ is a morphism of groups over $F$, $\alg{\phi}: \alg{G}_1' \rightarrow \alg{G}_2'$ making the following diagram commute:
$$\xymatrix{
1 \ar[r] & \alg{G}_{\mult} \ar[r]^{\alg{\iota}_1} \ar[d]^{=} & \alg{G}_1' \ar[r]^{\alg{p}_1} \ar[d]^{\alg{\phi}} & \alg{G} \ar[r] \ar[d]^{=} & 1 \\
1 \ar[r] & \alg{G}_{\mult} \ar[r]^{\alg{\iota}_2} & \alg{G}_2' \ar[r]^{\alg{p}_2} & \alg{G} \ar[r] & 1
}$$
This defines a category $\Cat{CExt}(\alg{G}, \alg{G}_{\mult})$ of central extensions of $\alg{G}$ by $\alg{G}_{\mult}$.  (Occasionally we might write $\Cat{CExt}_F(\alg{G}, \alg{G}_{\mult})$ to acknowledge the field of definition).
\end{definition}

\begin{proposition}
The category $\Cat{CExt}(\alg{G}, \alg{G}_{\mult})$ is a {\em groupoid}; every morphism in this category is an isomorphism.  The automorphism group of any object in this category is naturally isomorphic to the abelian group $X_F(\alg{G}) = \Hom_F(\alg{G}, \alg{G}_{\mult})$.
\end{proposition}
\proof
The fact that $\Cat{CExt}(\alg{G}, \alg{G}_{\mult})$ is a groupoid follows from a quick diagram chase.  As in Brylisnki-Deligne \cite{B-D}, and following Grothendieck \cite{SGA7}, the category of central extenions of $\alg{G}$ by $\alg{G}_{\mult}$ is equivalent to the category of multiplicative $\alg{G}_{\mult}$ torsors on $\alg{G}$.  The automorphisms of such torsors are the global multiplicative sections of the sheaf $(U \mapsto \alg{G}_{m,U})$ (for $U$ Zariski open in $\alg{G}$) over $\alg{G}$, i.e., regular functions from $\alg{G}$ to $\alg{G}_{\mult}$ which are multiplicative, i.e., the elements of $\Hom_F(\alg{G}, \alg{G}_{\mult})$.
\qed

Central extensions of $\alg{G}$ by $\alg{G}_{\mult}$ can be ``added'' via the Baer sum.
\begin{definition}
Let $(\alg{G}_1', \alg{p}_1, \alg{\iota}_1)$ and $(\alg{G}_2', \alg{p}_2, \alg{\iota}_2)$ be two central extensions of $\alg{G}$ by $\alg{G}_{\mult}$.  Let $\alg{\Delta}$ and $\alg{\nabla}$ denote the diagonal and antidiagonal embeddings of $\alg{G}_{\mult}$ into the center of the fibre product $\alg{G}_1' \times_{\alg{G}} \alg{G}_2'$.  The Baer sum $\alg{G}' = \alg{G}_1' + \alg{G}_2'$ is the quotient group:
$$\alg{G}' = \frac{\alg{G}_1' \times_{\alg{G}} \alg{G}_2'}{\alg{\nabla}(\alg{G}_{\mult})}.$$
The Baer sum is naturally a central extension of $\alg{G}$ by $\alg{G}_{\mult}$, with projection $\alg{p} = \alg{p}_1 + \alg{p}_2$ given by $\alg{p}_1$ on the first factor or equivalently $\alg{p}_2$ on the second factor, and with inclusion $\alg{\iota} = \alg{\iota}_1 + \alg{\iota}_2$ given by the diagonal embedding $\alg{\Delta}$.
\end{definition}
We refer to SGAIII, Expo.22, Section 4.3 \cite{SGA3III} for more on quotients of reductive groups by central tori, as used in the above construction.  The Baer sum $\alg{G}'$ is a reductive group over $F$ whose $L$-points (for a field $L$ containing $F$) are given by:
$$\alg{G}'(L) = [\alg{G}_1' + \alg{G}_2'](L) = \frac{  \{ (g_1', g_2') \in \alg{G}_1'(L) \times \alg{G}_2'(L) : g_1 = g_2 \} }{ \{ (z,z^{-1}) : z \in L^\times \} }.$$

This sum (defined above on objects of $\Cat{CExt}(\alg{G},\alg{G}_{\mult})$) extends to a functor:
$$+:  \Cat{CExt}(\alg{G},\alg{G}_{\mult}) \times \Cat{CExt}(\alg{G},\alg{G}_{\mult}) \rightarrow \Cat{CExt}(\alg{G},\alg{G}_{\mult}).$$
There are natural isomorphisms of functors which express the commutativity and associativity of the Baer sum.  A thorough way of describing the resulting structure on $\Cat{CExt}(\alg{G},\alg{G}_{\mult})$ is:
\begin{proposition}
The category $\Cat{CExt}(\alg{G},\alg{G}_{\mult})$, endowed with the Baer sum and natural commutativity and associativity isomorphisms, is a strictly commutative Picard groupoid (see Deligne, SGA IV \cite{SGA4}).
\end{proposition}

\subsection{Compatibilities}

Fix a central extension $(\alg{G}', \alg{p}, \alg{\iota})$ of $\alg{G}$ by $\alg{G}_{\mult}$.  The results of Section 4.3 of SGA III, Expo.22 \cite{SGA3III}, quickly imply the following:
\begin{proposition}
\label{ToriBijection}
If $\alg{T}$ is a maximal $F$-torus in $\alg{G}$, then its preimage $\alg{T}' = \alg{p}^{-1}(\alg{T})$ is a maximal torus in $\alg{G}'$.  This determines a bijection between the maximal $F$-tori in $\alg{G}$ and the maximal $F$-tori in $\alg{G}'$.
\end{proposition}

A crucial structural property of these central extensions is that they uniquely split over smooth unipotent subgroups:
\begin{theorem}
Let $\alg{U}$ be a smooth unipotent subgroup of $\alg{G}$ over $F$.  Then there exists a unique morphism of groups over $F$,  $\alg{s}: \alg{U} \rightarrow \alg{G}'$, such that $\alg{p} \circ \alg{s} = \alg{\Id}_{\alg{U}}$.  This morphism embeds $\alg{U}$ as a closed subgroup of $\alg{G}'$.
\label{UnipotentlyTrivialGm}
\end{theorem}
\proof
This directly follows from SGA III, Expo.17, Theorem 6.1.1 \cite{SGA3II}, since we assume $F$ is perfect.
\qed

This theorem is applicable to the most important examples of unipotent subgroups:
\begin{proposition}
Let $\alg{P_1}$ and $\alg{P_2}$ be parabolic subgroups of $\alg{G}$ over $F$.  Let $\alg{U_1}$ and $\alg{U_2}$ be the unipotent radicals of $\alg{P_1}$ and $\alg{P_2}$, respectively.  Then $\alg{U_1}$, $\alg{U_2}$, and $\alg{U_1} \cap \alg{U_2}$ are smooth unipotent subgroups of $\alg{G}$.
\end{proposition}
\proof
The smoothness of $\alg{U_1}$ and $\alg{U_2}$ follows from SGA III, Expo.26, Proposition 2.1 \cite{SGA3III}.  It is proven by identifying each of these unipotent groups (as a variety) with a product of smooth closed subgroups -- root subgroups -- of $\alg{G}$.  The two parabolic subgroups $\alg{P_1}$ and $\alg{P_2}$ contain a common maximal torus $\alg{T}$ (SGA III, Expo. 26, Lemma 4.1.1 \cite{SGA3III}); the intersection $\alg{U_1} \cap \alg{U_2}$ is a product of smooth closed root subgroups, with respect to this common torus.
\qed

Since $\alg{G}'$ is a central extension of $\alg{G}$, it follows that the conjugation action of $\alg{G}'$ on itself factors uniquely through the quotient $\alg{G}$:
$$\alg{\Int}: \alg{G} \times \alg{G}' \rightarrow \alg{G}'.$$
At the level of points, we write $[\Int(g)](x) = g' x (g')^{-1}$, where $g'$ is any lift of $g$.  In this way, $\alg{G}$ acts by conjugation on the variety of maximal tori in $\alg{G}'$, the variety of Borel subgroups in $\alg{G}'$, etc..

\begin{proposition}
The projection map $\alg{p}$ yields $\alg{G}$-equivariant isomorphisms over $F$ from:
\begin{enumerate}
\item
the variety $\alg{\Tor}(\alg{G}')$ of maximal tori in $\alg{G}'$ to the variety $\alg{\Tor}(\alg{G})$ of maximal tori in $\alg{G}$.
\item
the variety $\alg{\Bor}(\alg{G}')$ of Borel subgroups in $\alg{G}'$ to the variety $\alg{Bor}(\alg{G})$ of Borel subgroups in $\alg{G}$.
\item
the variety $\alg{\BorTor}(\alg{G}')$ of pairs $(\alg{B'}, \alg{T}')$ consisting of a Borel subgroup $\alg{B}'$ in $\alg{G}'$ and a maximal torus $\alg{T}'$ contained in $\alg{B}'$ to the corresponding variety $\alg{\BorTor}(\alg{G})$ of pairs in $\alg{G}$.
\item
the Springer variety $\alg{\Spr}(\alg{G}')$ of pairs $(\alg{B}', u)$ consisting of a Borel subgroup $\alg{B}'$ in $\alg{G}'$ and an element of its unipotent radical, to the corresponding Springer variety $\alg{\Spr}(\alg{G})$ of pairs in $\alg{G}$.
\end{enumerate}
\end{proposition}
\proof
Proposition \ref{ToriBijection} implies (1).  The varieties of Borel subgroups can be identified, as $\alg{G}$-varieties over $F$, with $\alg{G} / \alg{B}$ and $\alg{G}' / \alg{B}'$ (after a choice of base point), which are isomorphic via $\alg{p}$.  This demonstrates (2), and (3) is similar.  Theorem \ref{UnipotentlyTrivialGm} (or a version thereof, valid over a more general base variety) and (2) leads to a proof of (4).
\qed

The map $\alg{p}$ induces an isomorphism of Weyl groups, in every way possible:  first, if $\alg{T}$ is a maximal torus in $\alg{G}$, and $\alg{T'} = \alg{p}^{-1}(\alg{T})$, then $\alg{p}$ induces an isomorphism of finite \'etale groups over $F$:
$$\alg{W}(\alg{G}', \alg{T}') = \alg{\Norm}_{\alg{G}'}(\alg{T}') / \alg{T}' \rightarrow \alg{W}(\alg{G}, \alg{T}) = \alg{\Norm}_\alg{G}(\alg{T}) / \alg{T}.$$
This isomorphism is compatible with conjugation of tori, leading to an isomorphism from ``{\em the}'' Weyl group $\alg{W}$ of $\alg{G}'$ to ``{\em the}'' Weyl group of $\alg{G}$, in the sense of Section 1.1 of Deligne-Lusztig \cite{DeL}.

Since $\alg{p}$ induces a $\alg{G}$-equivariant isomorphism of varieties from $\alg{\Bor}(\alg{G}')$ to $\alg{\Bor}(\alg{G})$, it also induces a bijection from the $\alg{G}$-orbits on $\alg{\Bor}(\alg{G}') \times \alg{\Bor}(\alg{G}')$ to the $\alg{G}$-orbits on $\alg{\Bor}(\alg{G}) \times \alg{\Bor}(\alg{G})$.  In this way, $\alg{p}$ induces a bijection of Weyl groups, compatible with the Bruhat decomposition.  In particular, the bijection between the Borel subgroups of $\alg{G}'$ and those of $\alg{G}$ preserves the relation of ``being in relative position $w$'' for any $w$ in the Weyl group.

\subsection{Classification}

Let $F^{\sep}$ denote a separable closure of $F$ (which is an algebraic closure, since $F$ is perfect), and let $\Gamma = \Gal(F^{\sep} / F)$.  Let $\alg{T}$ be a maximal torus in $\alg{G}$, defined over $F$.  Let $(X,\Phi,Y,\Phi^\vee)$ denote the resulting (absolute) root system.  Thus $X$ and $Y$ are naturally $\ZZ[\Gamma]$-modules.  For a central extension $(\alg{G}', \alg{p}, \alg{\iota})$ of $\alg{G}$ by $\alg{G}_{\mult}$ as before, let $\alg{T}' = \alg{p}^{-1}(\alg{T})$, and $Y' = Y(\alg{T}')$.  This gives an extension of $\ZZ[\Gamma]$-modules that depends functorially on the central extension $(\alg{G}', \alg{p}, \alg{\iota})$:
\begin{equation}
0 \rightarrow \ZZ \rightarrow Y' \rightarrow Y \rightarrow 0.
\label{CentExtYZ}
\end{equation}
Somewhat more precisely,
\begin{proposition}
The above construction yields an additive functor (see \cite{SGA4}, Definition 1.4.5) of strictly commutative Picard groupoids:
$$\Cat{CExt}(\alg{G},\alg{G}_{\mult}) \rightarrow \Cat{Ext}_\Gamma(Y, \ZZ).$$
Here, $\Cat{Ext}_\Gamma(Y,\ZZ)$ denotes the category of extensions of $Y$ by $\ZZ$ in the abelian category of $\ZZ[\Gamma]$-modules.
\end{proposition}
\proof
The functoriality of this construction is clear.  Furthermore, the cocharacter lattice of the Baer sum is precisely the Baer sum of the cocharacter lattices, so this functor respects the Picard category structure.
\qed

We may refine this functor to obtain an equivalence of Picard groupoids.  The following theorem is analogous to the Main Theorem of Brylinski-Deligne \cite{B-D}.  The following theorem classifies central extensions of $\alg{G}$ by $\alg{G}_{\mult}$, while Brylinski and Deligne classify central extensions of $\alg{G}$ by $\alg{K}_2$.  The classifications are very similar in spirit, but the result below does not follow from \cite{B-D}, and the proof is quite different (and easier in our case).
\begin{theorem}
\label{CEByGm}
The category of central extensions of $\alg{G}$ by $\alg{G}_{\mult}$ is equivalent to the category of quadruples $(Y', p, \iota, \phi)$ as follows:  $(Y',p, \iota)$ is a $\ZZ[\Gamma]$-module extension of $Y$ by $\ZZ$:
$$\xymatrix{0 \ar[r] & \ZZ \ar[r]^\iota & Y' \ar[r]^{p} & Y \ar[r] & 0.}$$
Let $\alg{f}: \alg{G}_{\sconn} \rightarrow \alg{G}$ be the simply-connected cover of the derived group of $\alg{G}$, $\alg{T}_{\sconn} = \alg{f}^{-1}(\alg{T})$, and $Y_{\sconn}$ the cocharacter group of $\alg{T}_{\sconn}$.  The last part of the quadruple, $\phi$, is a Galois-equivariant morphism from $Y_{\sconn} \times \ZZ$ to $Y_{\sconn}$, making the following diagram commute:
$$\xymatrix{
0 \ar[r] & \ZZ \ar[r] \ar[d]^= & Y_{\sconn} \times \ZZ \ar[r] \ar[d]^{\phi} & Y_{\sconn} \ar[r]  \ar[d]^{f_\ast} & 0 \\
0 \ar[r] & \ZZ \ar[r]^{\iota} & Y' \ar[r]^p & Y \ar[r] & 0.
}$$
Morphisms from a quadruple $(Y_1', p_1, \iota_1, \phi_1)$ to a quadruple $(Y_2', p_2, \iota_2, \phi_2)$ are morphisms of $\ZZ[\Gamma]$-modules from $Y_1'$ to $Y_2'$ making the large but obvious diagram of $\ZZ[\Gamma]$-modules commute.
\end{theorem}
\proof
If $(\alg{G}', \alg{p}, \alg{\iota})$ is a central extension of $\alg{G}$ by $\alg{G}_{\mult}$, then the cocharacter lattices yield an extension of $\ZZ[\Gamma]$-modules $(Y', p, \iota)$ as above.  Furthermore, the pullback of the central extension yields a central extension $\alg{G}_{\sconn}'$ of $\alg{G}_{\sconn}$ by $\alg{G}_{\mult}$, which splits uniquely (since $\alg{G}_{\sconn}$ is simply-connected).  Letting $Y_{\sconn}'$ be the cocharacter lattice of the maximal torus $\alg{T}_{\sconn}'$ in $\alg{G}_{\sconn}'$ (the pullback of $\alg{T}'$, from $\alg{G}'$ to $\alg{G}_{\sconn}'$), we find that $Y_{\sconn}'$ is canonically identified with $Y_{\sconn} \times \ZZ$.  The covering map from $\alg{G}_{\sconn}'$ to $\alg{G}'$ yields the requisite map $\phi$ from $Y_{\sconn} \times \ZZ$ to $Y'$.

This describes the functor from $\Cat{CExt}(\alg{G}, \alg{G}_{\mult})$ to the category of quadruples.  It is compatible with the Baer sum as well.  To prove that this functor is an equivalence, we prove first that it is bijective on automorphism groups; this implies that the functor is fully faithful, since both categories are groupoids.

The automorphism group of a central extension $(\alg{G}', \alg{p}, \alg{\iota})$ of $\alg{G}$ by $\alg{G}_{\mult}$ can be identified with $\Hom_F(\alg{G}, \alg{G}_{\mult})$.  This group of $F$-rational characters of $\alg{G}$ embeds (by restriction, naturally) as a subgroup of $X_F(\alg{T}) = \Hom_F(\alg{T}, \alg{G}_{\mult})$.  The image of this embedding is given by:
$$\Hom_F(\alg{G}, \alg{G}_{\mult}) \isom \Ker(X_F(\alg{T}) \rightarrow X_F(\alg{T}_{\sconn})) \isom \Hom_{\Gamma} (Y / f_\ast Y_{\sconn}, \ZZ).$$
On the other hand, the automorphisms of a quadruple $(Y', p, \iota, \phi)$ are precisely the automorphisms of an extension $\ZZ \rightarrow Y' \rightarrow Y$ of $\ZZ[\Gamma]$-modules which pull back to the trivial automorphism of an extension $\ZZ \rightarrow Y_{\sconn} \times \ZZ \rightarrow Y_{\sconn}$.  Such automorphisms are naturally identified with elements of $\Hom_{\Gamma}(Y, \ZZ)$ which pull back to trivial elements of $\Hom_{\Gamma}(Y_{\sconn}, \ZZ)$.  Hence this automorphism group is naturally identified with $\Hom_{\Gamma}(Y / f_\ast Y_{\sconn}, \ZZ)$. Hence our functor is bijective on automorphism groups (leaving the reader to check that a diagram of isomorphisms commutes).

Now to prove essential surjectivity of this functor, we may assume $\alg{G}$ is split by \'etale descent, since we have verified compatibility with automorphism groups.  In Section 2.4 of \cite{Kot}, Kottwitz demonstrates an isomorphism, functorial for ``normal'' (Section 1.8 of \cite{Kot}) homomorphisms,
$$\Pic(\alg{G}) \isom \pi_0 Z(\alg{G}^\vee).$$
By Hilbert's Theorem 90, line bundles on $\alg{G}$ can be rigidified at the identity element.  Since the projection and multiplication maps 
$$\pr_1, \pr_2, m: \alg{G} \times \alg{G} \rightarrow \alg{G}$$
are normal in the sense of \cite{Kot}, the Kottwitz isomorphisms are compatible:
$$\xymatrix{
\Pic(\alg{G}) \ar[d]^{\pr_1^\ast, \pr_2^\ast, m^\ast} \ar[rr]^{\Kott} & &  \pi_0 Z(\alg{G}^\vee) \ar[d]^{\pr_1^\vee, \pr_2^\vee, m^\vee} \\
\Pic(\alg{G} \times \alg{G}) \ar[rr]^{\Kott} & & \pi_0 Z(\alg{G}^\vee \times \alg{G}^\vee).
}$$
On the right side, we find easily that $\pr_1^\vee(z) \cdot \pr_2^\vee(z) = m^\vee(z)$; on the left side, therefore, we find the same equality in Picard groups; $\pr_1^\ast(L) \cdot \pr_2^\ast(L) = m^\ast(L)$, for an invertible sheaf $L$ on $\alg{G}$.  It follows from Proposition 4.2 of SGA 7, Expo.VIII \cite{SGA7}, that the line bundles classified by $\Pic(\alg{G})$ define extensions of $\alg{G}$ by $\alg{G}_{\mult}$.  In other words, $\Pic(\alg{G})$ is naturally isomorphic to the group $\CExt(\alg{G}, \alg{G}_{\mult})$ of isomorphism classes in $\Cat{CExt}(\alg{G}, \alg{G}_{\mult})$.  

Thus it remains to prove that this group of isomorphism classes $\Pic(\alg{G})$ -- naturally isomorphic to $\pi_0 Z(\alg{G}^\vee)$ on one hand -- is isomorphic to the group of isomorphism classes of quadruples $(Y', p, \iota, \phi)$ discussed above.  The isomorphism classes of quadruples $(Y', p, \iota, \phi)$ are classified by the hypercohomology of the two-term complex $f_\ast: Y_{\sconn} \rightarrow Y$, with coefficients in $\ZZ$:
$$\HH^2(Y_{\sconn} \rightarrow Y, \ZZ) \isom \Ext^1(Y / f_\ast Y_{\sconn}, \ZZ).$$
(Compare to (6.3.1) of \cite{B-D}).  This is isomorphic, by Proposition 1.10 of \cite{Bor} and Lemma 2.2 of \cite{Kot}, to $\pi_0 Z(\alg{G}^\vee)$ as required.  Again we leave it to the reader to verify that this isomorphism agrees with the one given by our functor.  This is not as much of a ``cop out'' as it might seem -- the maps occurring in the Kottwitz isomorphism, and in the work of Borovoi, are also given by considering maps of cocharacter lattices, and so agreement is inevitable.
\qed

\begin{remark}
The identification of $\Pic(\alg{G})$ with $\CExt(\alg{G}, \alg{G}_{\mult})$ is also proven, without recourse to the dual group, by Colliot-Th\`el\'ene in Theorem 5.6 of \cite{C-T}.  This statement was almost certainly known decades ago to the experts; there are similarities to Chapter VII of Raynaud's thesis \cite{Ray}.  It also appears in an unpublished communication of O. Gabber.  We thank Mikhail Borovoi and Brian Conrad for providing these references.
\end{remark}

This theorem describes, up to equivalence, the category of central extensions of $\alg{G}$ by $\alg{G}_{\mult}$.  Such a description is useful for purposes of descent, and for tracing how an automorphism of central extensions (as would be induced by changing a cocycle within a cohomology class) affects other parameters.

\begin{corollary}
If $\alg{G}$ is a semisimple group over $F$, then every central extension of $\alg{G}$ by $\alg{G}_{\mult}$ is rigid, i.e., has no nontrivial automorphisms.
\end{corollary}
\proof
Let $(\alg{G}', \alg{p}, \alg{\iota})$ be a central extension of $\alg{G}$ by $\alg{G}_{\mult}$.  Its automorphism group is isomorphic to $\Hom_F(\alg{G}, \alg{G}_{\mult})$, which is trivial when $\alg{G}$ is semisimple (recall we always assume $\alg{G}$ to be connected).
\qed

\begin{example}
Let $\alg{G} = \alg{PGL}_2$.  Then all central extensions of $\alg{G}$ are rigid.  The isomorphism classes of such central extensions are in natural bijection with $\pi_0 Z(\alg{SL}_2) \isom \mu_2$.  These two isomorphism classes of central extensions are represented by the two familiar extensions:
$$1 \rightarrow \alg{G}_{\mult} \rightarrow \alg{G}_{\mult} \times \alg{PGL}_2 \rightarrow \alg{PGL}_2 \rightarrow 1,$$
$$1 \rightarrow \alg{G}_{\mult} \rightarrow \alg{GL}_2 \rightarrow \alg{PGL}_2 \rightarrow 1.$$
\end{example}

\section{Unipotently split extensions}

In this section, we consider another kind of central extension which is a convenient compromise between ``abstract'' (in the terminology of \cite{Ti1}) group theory and algebraic group theory.  This compromise avoids the potential trouble of having too many central extensions in abstract group theory, and avoids the hard work (as in \cite{Moo}) in classifying central extensions in a topological category.  This compromise -- our class of unipotently split extensions -- arises naturally from the central extensions of reductive groups by $\alg{K}_2$ -- the class of extensions studied by Brylinski and Deligne \cite{B-D}.

In this section, we allow $F$ to be any field (perfect or not).  We briefly allow $\alg{G}$ to be any algebraic group over $F$, and $G = \alg{G}(F)$ the group of points.  We fix an abelian group $\mu$, and study a class of central extensions $(\tilde G, p, \iota)$ of $G$ by $\mu$:
$$\xymatrix{
1 \ar[r] & \mu \ar[r]^{\iota} & {\tilde G} \ar[r]^{p} & G \ar[r] &1}.$$
In this section, we consider such central extensions endowed with a ``unipotent splitting'':
\begin{definition}
\label{UnipotentSplitting}
Let $(\tilde G, p, \iota)$ be a central extension of $G$ by $\mu$.  A {\em unipotent splitting} of $(\tilde G, p, \iota)$ is a family of homomorphisms $\{ \tilde \eta: U \rightarrow \tilde G \}$ indexed by all homomorphisms $\alg{\eta}: \alg{U} \rightarrow \alg{G}$ from split unipotent groups to $\alg{G}$, defined over $F$, satisfying the following conditions:
\begin{enumerate}
\item
For each $\alg{\eta}: \alg{U} \rightarrow \alg{G}$, $p \circ \tilde \eta = \eta$ as homomorphisms from $U$ to $G$:
$$\xymatrix{U \ar@/_/[rr]_{\eta} \ar[r]^{\tilde \eta} & \tilde G \ar[r]^p & G.}$$
\item
For every pair $\alg{U}_1$, $\alg{U}_2$ of split unipotent groups, and commutative diagram of groups over $F$:
$$\xymatrix{
\alg{U_1} \ar@/_/[rr]_{\alg{\eta}_1} \ar[r]^{\alg{f}} & \alg{U}_2 \ar[r]^{\alg{\eta}_2} & \alg{G},}$$
the homomorphisms $\tilde \eta_2$ and $\tilde \eta_1$ satisfy $\tilde \eta_2 \circ f = \tilde \eta_1$:
$$\xymatrix{
U_1 \ar@/_/[rr]_{\tilde \eta_1} \ar[r]^f & U_2 \ar[r]^{\tilde \eta_2} & \tilde G,}$$
\item
For each homomorphism from a split unipotent group, $\alg{\eta}: \alg{U} \hookrightarrow \alg{G}$, and each element $g \in G$, the following diagram commutes:
$$\xymatrix{
U \ar[r]^{\tilde \eta} \ar[dr]_{ [\Int(g) \circ \eta]^{\sim}} & \tilde G \ar[d]^{\Int(g)} \\
& \tilde G
}.$$
\end{enumerate}
\end{definition}
\begin{remark}
There may be a set-theoretic subtlety in defining a ``family'' of homomorphisms indexed by the ``set'' of homomorphisms from split unipotent groups into $\alg{G}$.  This is easily resolved by restricting to a sufficient set of split unipotent groups.  
\end{remark}

\subsection{Chevalley groups}

Assume for now that $\alg{G}$ is a split, semisimple, simply-connected group over $F$, and $\alg{S}$ is an $F$-split maximal torus in $\alg{G}$.  Let $\Phi = \Phi(\alg{G}, \alg{S})$ denote the resulting set of roots, and $(X, \Phi, Y, \Phi^\vee)$ the root datum.  For $\alpha \in \Phi$, the associated reflections of $X$ and $Y$ are defined by:
$$s_\alpha(x) = x - \langle \alpha^\vee, x \rangle \alpha, \quad s_{\alpha^\vee}(y) = y - \langle y, \alpha \rangle \alpha^\vee,$$
for all $x \in X, y \in Y$.  

Since $\alg{G}$ is simply-connected, the cocharacter lattice $Y$ is generated as a $\ZZ$-module by $\Phi^\vee$.  Even more, Brylinski and Deligne prove (Lemma 11.5 of \cite{B-D}) that $Y$ can be presented as the quotient of the free abelian group $\ZZ \langle \alpha^\vee \rangle_{\alpha \in \Phi}$ modulo the relations arising from root reflections:
$$s_\alpha(\beta)^\vee = \beta^\vee - \langle \beta^\vee, \alpha \rangle \alpha^\vee.$$

By SGA 3, Expo. 23, Proposition 6.2 \cite{SGA3III}, we may choose a {\em Chevalley system} (see Definition 6.1, loc.~cit.) on $\alg{G}$.  Such a system yields a set $\{ \alg{e}_\alpha: \alg{G}_{\add} \rightarrow \alg{U}_\alpha \}$ of isomorphisms from the additive group $\alg{G}_{\add}$ onto the $\alpha$ root subgroup $\alg{U}_\alpha$ for each root $\alpha \in \Phi$.  Define a map (not a homomorphism) $\alg{n}_\alpha: \alg{G}_{\mult} \rightarrow \alg{G}$, for $\alpha \in \Phi$, by
$$n_\alpha(z) = e_\alpha(z) e_{-\alpha}(-z^{-1}) e_\alpha(z).$$
From the definition of Chevalley system, the following identity holds for all $\alpha, \beta \in \Phi$, and all $u \in F$:
$$n_\alpha(1) e_\beta(u) n_\alpha(1)^{-1} = e_{s_\alpha(\beta)}(\pm u),$$
where the sign depends only on $\alpha$ and $\beta$.  

Choose also a system of positive roots, yielding a partition $\Phi = \Phi^+ \cup \Phi^-$ and a set of simple roots $\Delta \subset \Phi^+$.  For each positive root $\alpha$, the Chevalley system yields a central isogeny $\alg{\phi}_\alpha: \alg{SL}_2 \rightarrow \alg{G}_\alpha$, where $\alg{G}_\alpha$ is a closed subgroup of $\alg{G}$ containing $\alg{U}_\alpha$ and $\alg{U}_{-\alpha}$, and for which
$$\phi_\alpha \Matrix{1}{u}{0}{1} = e_\alpha(u).$$
Since it is convenient, and we have our choice of signs, we require the Chevalley system to satisfy the identity:
$$\phi_{\alpha} \Matrix{1}{0}{u}{1} = e_{-\alpha}(u), \mbox{ for all } \alpha \in \Phi^+, u \in F.$$
It follows that
$$\phi_\alpha \Matrix{0}{z}{-z^{-1}}{0} = n_\alpha(z), \mbox{ for all } z \in F^\times.$$

The following ``Chevalley-Steinberg'' relations hold:
\begin{proposition}
Relations (B) and (B') from Section 6 of \cite{Ste}, hold in the group $G$:
\begin{enumerate}
\item[(B)]
For all roots $\alpha, \beta \in \Phi$, such that $\alpha \neq \pm \beta$, there is an ordering of the set of roots of the form $\{ i \alpha + j \beta \}_{0 < i,j \in \ZZ}$, and integers $c_{ij}(\alpha, \beta)$, such that:
$$[e_\alpha(u), e_\beta(v)] = \prod_{i \alpha + j \beta} e_{i \alpha + j \beta}(c_{ij}(\alpha, \beta) u^i v^j ).$$
\item[(B')]
For all roots $\alpha \in \Phi$, all $z \in F^\times$, $u \in F$,
$$[\Int(n_\alpha(z))] \left( e_\alpha(u) \right) = n_\alpha(z) e_\alpha(u) n_\alpha(-z) = e_{-\alpha}(-z^{-2} u).$$
\end{enumerate}
\end{proposition}
\proof
For the relation (B) (and more precise information about it), we refer to Proposition 6.4 of SGA 3, Expo.23 \cite{SGA3III}.  The relation (B') follows directly from our sign convention, and the corresponding relation in $\alg{SL}_2(F)$. 
\qed

Now we consider a central extension $(\tilde G, p, \iota)$ of $G$ by an abelian group $\mu$.  We moreover suppose that this central extension is endowed with a unipotent splitting $\{ \tilde \eta: U \rightarrow \tilde G \}$.  The  homomorphisms $\alg{e}_\alpha: \alg{G}_{\add} \rightarrow \alg{G}$ then lift to homomorphisms $\tilde e_\alpha: F \rightarrow \tilde G$.  Define, for $z \in F^\times$,
$$\tilde n_\alpha(z) =  \tilde e_\alpha(z) \tilde e_{-\alpha}(-z^{-1}) \tilde e_\alpha(z).$$
\begin{theorem}
\label{ChevalleySteinberg}
The Chevalley-Steinberg relations of the previous proposition hold, with $G$ replaced by $\tilde G$, $e_\alpha$ replaced by $\tilde e_\alpha$ and $n_\alpha$ replaced by $\tilde n_\alpha$.
\end{theorem}
\proof
For the relation (B), consider the split unipotent subgroup $\alg{U} \subset \alg{G}$ spanned by the root spaces $\alg{U}_{i \alpha + j \beta}$ for non-negative integers $i,j$.  The maps $\alg{e}_{i \alpha + j \beta}: \alg{G}_{\add} \rightarrow \alg{G}$ factor through $\alg{U}$, and we define $\alg{f}_{i \alpha + j \beta}: \alg{G}_{\add} \rightarrow \alg{U}$ and $\alg{\eta}: \alg{U} \hookrightarrow \alg{G}$ so that $\alg{e}_{i \alpha + j \beta} = \alg{\eta} \circ \alg{f}_{i \alpha + j \beta}$.  It follows from (2) in Definition \ref{UnipotentSplitting} that for all $u,v \in F$,
\begin{eqnarray*}
[\tilde e_\alpha(u), \tilde e_\beta(v)] & = & [\tilde \eta \circ f_\alpha(u), \tilde \eta \circ f_\beta(v)] \\
& = & \tilde \eta \left( [f_\alpha(u), f_\beta(v)] \right) \\
& = & \tilde \eta \left( \prod_{i \alpha + j \beta} f_{i \alpha + j \beta}(c_{ij} u^i v^j ) \right) \\
& = & \prod_{i \alpha + j \beta} \left( \tilde \eta \circ f_{i \alpha + j \beta}(c_{ij} u^i v^j ) \right) \\
& = & \prod_{i \alpha + j \beta} \tilde e_{i \alpha + j \beta}(c_{ij} u^i v^j ).
\end{eqnarray*}

By relation (B') in the group $G$, the following diagram of groups and homomorphisms over $F$ commutes (for any $z \in F^\times$):
$$\xymatrix{
\alg{G}_{\add} \ar[r]^{\alg{e}_\alpha} \ar[d]_{\alg{\Mult}(-z^{-2})} & \alg{G} \ar[d]^{\alg{\Int}(n_\alpha(z))} \\
\alg{G}_{\add} \ar[r]^{\alg{e}_{-\alpha}} & \alg{G}.
}$$
It follows from (2) and (3) in Definition \ref{UnipotentSplitting} that:
\begin{eqnarray*}
[\Int(n_\alpha(z))] \circ \tilde e_\alpha & = & \left( \Int(n_\alpha(z)) \circ e_\alpha \right)^{\sim} \\
& = & \left( e_{-\alpha} \circ \Mult(-z^{-2}) \right)^{\sim} \\
& = & \tilde e_{-\alpha} \circ \Mult(-z^{-2}).
\end{eqnarray*}
This demonstrates that relation (B') holds in $\tilde G$ as well.
\qed

We carry on to derive relations in the central extension $\tilde G$.  The central extension $\tilde G$ restricts to a central extension of the $F$-points of the split torus $\alg{S}$:
$$1 \rightarrow \mu \rightarrow \tilde S \rightarrow S \rightarrow 1.$$
Define elements of $\tilde S$ by 
$$\tilde h_\alpha(z) = \tilde n_\alpha(z) \tilde n_\alpha(-1), \mbox{ for all } \alpha \in \Phi, z \in F^\times.$$
These elements project onto $h_\alpha(z) = \alpha^\vee(z) \in S$, for every root $\alpha \in \Phi$.  The maps $\tilde h_\alpha$ are not necessarily homomorphisms, and the deviation is measured by a 2-cocycle for every root:
$$\sigma_\alpha(z_1, z_2) = \tilde h_\alpha(z_1) \tilde h_\alpha(z_2) \tilde h_\alpha(z_1 z_2)^{-1}.$$
Note that $\sigma_\alpha \in Z^2(F^\times, \mu)$, and $\sigma_\alpha(1, z) = \sigma_\alpha(z,1) = 1$ for all $z \in F^\times$.

\begin{corollary}
Choose two roots $\alpha, \beta \in \Phi$, and any elements $u \in F$, $z,v \in F^\times$.  Let $\gamma = s_\alpha(\beta) = \beta - \langle \alpha^\vee, \beta \rangle \alpha$.  Then for some constant $\epsilon = \epsilon(\alpha, \beta) = \pm 1$ (independent of $z,v$), the following relations hold:
\begin{eqnarray*}
\tilde n_\alpha(z) \tilde e_\beta(u) \tilde n_\alpha(z)^{-1} & = & \tilde e_{\gamma}(\epsilon z^{-\langle  \alpha^\vee, \beta \rangle} u), \\
\tilde n_\alpha(z) \tilde h_\beta(v) \tilde n_\alpha(z)^{-1} & = & \tilde h_\gamma(v) \sigma_\gamma(v,\epsilon z^{-\langle \alpha^\vee, \beta \rangle})^{-1}, \\
\tilde h_\alpha(z) \tilde e_\beta(u) \tilde h_\alpha(z)^{-1} & = & \tilde e_\beta ( z^{\langle \alpha^\vee, \beta \rangle} u ), \\
\tilde h_\alpha(z) \tilde h_\beta(v) \tilde h_\alpha(z)^{-1} & = & \tilde h_\beta( v ) \sigma_\beta(v, z^{\langle \alpha^\vee, \beta \rangle} )^{-1}.
\end{eqnarray*}
\end{corollary}
\proof
This follows immediately from our earlier Theorem \ref{ChevalleySteinberg}, and Lemma 37 in Steinberg's Yale lectures \cite{Ste}.
\qed

\subsection{Simply-connected semisimple groups}
Suppose that $\alg{G}$ is a simply-connected quasi-split semisimple group over $F$.  Then $\alg{G}$ is a direct product of simply-connected, quasi-split, $F$-almost-simple $F$-subgroups.  Each of these $F$-almost-simple factors is $F$-isomorphic to $\alg{R}_{L/F} \alg{H}$ (Weil restriction of scalars), for some finite separable extension $L/F$, and for some simply-connected, quasi-split, {\em absolutely} almost-simple group $\alg{H}$ over $L$.  This follows from \cite{BoT}, Section 6.21 (ii).

Since a central extension of a direct product of perfect groups is determined by a collection of central extensions of factor groups, we assume in this section that $\alg{G} = \alg{R}_{L/F} \alg{H}$, with $\alg{H}$ simply-connected, quasi-split, and absolutely almost-simple over $L$, as described above.  Thus $\alg{H}$ belongs to one of the following types:
$$\Type{A}_n, {}^2 \Type{A}_n, \Type{B}_n, \Type{C}_n, \Type{D}_n, {}^2 \Type{D}_n, {}^3 \Type{D}_4, {}^6 \Type{D}_4, \Type{E}_6, {}^2 \Type{E}_6, \Type{E}_7, \Type{E}_8, \Type{F}_4, \Type{G}_2.$$

Let $\alg{S}$ be a maximal $F$-split torus in $\alg{G}$, and $\Phi$ the resulting set of {\em relative} roots.  Let $\alg{W} = \alg{W}(\alg{G}, \alg{S})$ denote the resulting relative Weyl group.  Fix a system of positive roots, so that $\Phi = \Phi^+ \sqcup \Phi^-$.  Since we do not assume $\alg{G}$ to be split, it is possible that $\Phi$ is not reduced and the root spaces in $\Lie{g}$ may have dimension greater than $1$.  Define
$$\Phi_1 = \{ \alpha \in \Phi : \alpha/2 \not \in \Phi \}, \quad \Phi_2 = \{ \alpha \in \Phi : 2 \alpha \not \in \Phi \}.$$
Thus $\Phi_1$ is the set of {\em indivisible roots}, and $\Phi_2$ the set of {\em undoublable roots}.

For any root $\alpha$, we write $\alg{U}_\alpha$ for the unique unipotent subgroup of $\alg{G}$ whose Lie algebra is the sum of root spaces, for roots which are positive integer multiples of $\alpha$:
$$\Lie{u}_\alpha = \bigoplus_{k > 0 } \Lie{g}_{k \alpha}.$$
By Corollary 3.18 of \cite{BoT}, it is known that $\alg{U}_\alpha$ is a unipotent subgroup of $\alg{G}$, split over $F$.

When $\alpha \in \Phi_1^+$, let $\alg{G}_\alpha$ be the smallest closed subgroup of $\alg{G}$ containing $\alg{U}_\alpha$ and $\alg{U}_{-\alpha}$.  Then $\alg{G}_\alpha$ is a quasi-split,  almost-simple (over $F$) group, of $F$-rank $1$.  We consider two cases:
\subsubsection{When $\alpha \in \Phi_2$}  
When $\alpha \in \Phi$ and $2 \alpha \not \in \Phi$, there is a central isogeny over $F$:
$$\alg{\phi}_\alpha: \alg{R}_{E/F} \alg{SL}_{2,E} \rightarrow \alg{G}_\alpha,$$
where $E$ is a finite separable extension of $L$.    We choose this central isogeny so that the ``diagonal torus'' lands in $\alg{S}$, and an isomorphism $\alg{e}_\alpha: \alg{R}_{E/F} \alg{G}_{\add} \rightarrow \alg{U}_\alpha$ is given by:
$$\phi_\alpha \Matrix{1}{u}{0}{1} = e_\alpha(u), \mbox{ for all } u \in \alg{R}_{E/F} \alg{G}_{\add}(F) = E.$$
An isomorphism $\alg{e}_{-\alpha}: \alg{R}_{E/F} \alg{G}_{\add} \rightarrow \alg{U}_{-\alpha}$ is given by
$$\phi_\alpha \Matrix{1}{0}{u}{1} = e_{-\alpha}(u).$$
Define a map (not a homomorphism) of varieties over $F$, $\alg{n_\alpha}: \alg{R}_{E/F} \alg{G}_{\mult} \rightarrow \alg{G}$ by
$$n_\alpha(z) = e_\alpha(z) e_{-\alpha}(-z^{-1}) e_\alpha(z) = \phi_\alpha \Matrix{0}{z}{-z^{-1}}{0},$$
for all $z \in E^\times$.  Such elements represent the relative Weyl group reflection $s_\alpha$ in $\Norm_G(S)$.  A short calculuation gives:
\begin{equation}
\label{NErelation1}
n_\alpha(z) e_\alpha(u) n_\alpha(z)^{-1} = e_{-\alpha}(-z^{-2} u).
\end{equation}

Define a homomorphism $\alg{h}_\alpha: \alg{R}_{E/F} \alg{G}_{\mult} \rightarrow \alg{G}$ by $h_\alpha(z) = n_\alpha(z) n_\alpha(-1)$, for $z \in E^\times$.  Then
$$h_\alpha(z) = \phi_\alpha \Matrix{z}{0}{0}{z^{-1}}.$$
When restricted to $F^\times \subset E^\times$, $h_\alpha$ coincides with a cocharacter $\alpha^\vee$ of $\alg{S}$.  

\subsubsection{When $\alpha \not \in \Phi_2$}
When $\alpha \in \Phi^+$ and $2 \alpha$ is also a root, there is a central isogeny over $F$:
$$\alg{\phi}_\alpha: \alg{R}_{L/F} \alg{SU}_{3, E/L}  \rightarrow \alg{G}_\alpha,$$
where $E$ is a separable quadratic extension of $L$.  We write $\sigma$ for the nontrivial element of $\Gal(E/L)$.  Here the quasisplit group $\alg{SU}_{3, E/L}$ is given by:
$$\alg{SU}_{3,E/L}(L) = \{ g \in \alg{SL}_3(E) : g \Xi (g^\sigma)^t = \Xi \}, \mbox{ where } \Xi = \TMatrix{0}{0}{1}{0}{1}{0}{1}{0}{0}.$$
We chose $\alg{\phi}_\alpha$ so that the ``diagonal torus'' gets mapped into $\alg{S}$, and the upper-triangular (resp. lower-triangular) unipotent subgroups of $\alg{SU}_{3, E/L}$ get mapped to the unipotent subgroups $\alg{U}_\alpha$ (resp. $\alg{U}_{-\alpha}$).  Define a split unipotent group over $L$ by:
$$\alg{J}_{E/L}(L) = \{ (p, \ell) \in E^2 \mbox{ such that } \ell + \ell^\sigma + p p^\sigma  = 0 \},$$
where the group law is given by:
$$(p_1, \ell_1) \cdot (p_2, \ell_2) = (p_1 + p_2, \ell_1 + \ell_2 - p_1^\sigma p_2).$$
Of course, this only defines the $L$-points of this group, but it is easy to extend the above to the points over any $L$-algebra.

We define homomorphisms of groups over $F$, $\alg{e}_{\pm \alpha}: \alg{R}_{L/F} \alg{J}_{E/L} \rightarrow \alg{U}_{\pm \alpha}$ (following Deodhar \cite{Deo} but with choices made by Tits \cite{Tit}):
$$e_\alpha(p, \ell)  =  \phi_\alpha \TMatrix{1}{-p^\sigma}{\ell}{0}{1}{p}{0}{0}{1}, \quad e_{-\alpha}(p, \ell) =  \phi_\alpha \TMatrix{1}{0}{0}{p}{1}{0}{\ell}{-p^\sigma}{1}.$$
Let $\alg{J}_{E/L}^\ast = \alg{J}_{E/L} - \{ (0,0) \}$ denote the complement of the identity -- a subvariety of $\alg{J}_{E/L}$ defined over $L$.  Define a map of varieties over $F$, $\alg{n}_\alpha: \alg{R}_{L/F} \alg{J}_{E/L}^\ast \rightarrow \alg{G}$ by:
$$n_\alpha(c, d) = e_{-\alpha}(-c d^{-1}, d^{-\sigma}) e_\alpha(c,d) e_{-\alpha}(-c d^{-\sigma}, d^{-\sigma}).$$
Here, we note that $(0,0) \neq (c,d) \in E^2$ satisfies $c + c^\sigma + d d^\sigma = 0$, which implies that $d \neq 0$, and we write $d^{-\sigma}$ for $(d^{-1})^\sigma = (d^\sigma)^{-1}$.  A short computation demonstrates that
$$n_\alpha(c, d) = \phi_\alpha \TMatrix{0}{0}{d}{0}{-d^\sigma / d}{0}{d^{-\sigma}}{0}{0}.$$
Such elements represent the relative Weyl group reflection $s_\alpha$ in $\Norm_G(S)$.  

A computation gives, for all $(c,d) \in \alg{J}_{E/L}^\ast(L)$ and all $(p, \ell) \in \alg{J}_{E/L}(L)$:
\begin{eqnarray}
n_\alpha(c, d) e_\alpha(p,\ell) n_\alpha(c,d)^{-1} & = & e_{-\alpha} \left( -\frac{d^\sigma p}{ d^2}, \frac{\ell}{d d^\sigma} \right), \\
n_\alpha(c,d) e_{-\alpha}(p,\ell) n_\alpha(c,d)^{-1} & = & e_\alpha \left( \frac{d^2 p^\sigma}{d^\sigma}, \ell d d^\sigma \right).
\end{eqnarray}

We now follow the work of Deodhar, Section 2.11 of \cite{Deo}, for guidance.  If $\characteristic(L) \neq 2$, then there exists $\theta$ such that $E = L(\theta)$ and $\theta + \theta^\sigma = 0$.  Such a $\theta$ is unique up to scaling by $L^\times$, and we choose such a $\theta$ for what follows.  If $\characteristic(L) = 2$, then note that for every $\ell \in L$, $\ell + \ell^\sigma = 2 \ell = 0$.  Thus if $\characteristic(L) = 2$, we choose $\theta = 1$.

Define a homomorphism of groups over $F$, $\alg{e}_{2 \alpha}: \alg{R}_{L/F} \alg{G}_{\add} \rightarrow \alg{G}$ by
$$e_{2 \alpha}(\ell) = e_\alpha(0, \ell \theta) = \phi_\alpha \TMatrix{1}{0}{\ell \theta}{0}{1}{0}{0}{0}{1}.$$
Define $\alg{e}_{-2 \alpha}$ similarly by $e_{- 2 \alpha}(\ell) = e_{-\alpha}(0, \ell \theta)$.

Define a map of varieties over $F$, $\alg{n}_{2 \alpha}: \alg{R}_{L/F} \alg{G}_{\mult} \rightarrow \alg{G}$ by:
$$n_{2 \alpha}(\ell) = n_\alpha(0, \ell \theta).$$
Note that 
\begin{eqnarray*}
n_{2 \alpha}(\ell) & = & e_{-\alpha}(0, \ell^{-\sigma} \theta^{-\sigma}) e_\alpha(0, \ell \theta) e_{-\alpha}(0, \ell^{-\sigma} \theta^{-\sigma}) \\
& = & e_{-2 \alpha} \left( \frac{\ell^{-\sigma}}{\theta \theta^{\sigma}} \right) e_{2 \alpha}(\ell) e_{-2 \alpha} \left( \frac{\ell^{-\sigma}}{ \theta \theta^{\sigma}} \right).
\end{eqnarray*}
Finally, define an algebraic homomorphism $\alg{h}_{2 \alpha}: \alg{R}_{L/F} \alg{G}_{\mult} \rightarrow \alg{G}$ by:
$$h_{2 \alpha}(\ell) = n_{2 \alpha}(\ell) n_{2 \alpha}(-1).$$
A computation yields
$$h_{2 \alpha}(\ell) = \phi_\alpha \TMatrix{\ell}{0}{0}{0}{1}{0}{0}{0}{\ell^{-1}}.$$

\subsubsection{Covers}

Now, at last, we consider a  central extension $\tilde G$ of $G = \alg{G}(F)$ by $\mu$, endowed with a {\em unipotent splitting}, where $\alg{G} = \alg{R}_{L/F} \alg{H}$, and $\alg{H}$ is a quasisplit, simply-connected, semisimple, and absolutely almost simple group, defined over $L$.  Let $\alg{S}$ be a maximal $F$-split torus in $\alg{G}$, and $\alpha \in \Phi(\alg{G}, \alg{S})$ an indivisible (relative) root.  

Since $\alg{e}_\alpha$ is a homomorphism from a split unipotent group into $\alg{G}$, it lifts to a homomorphism:
$$\tilde e_\alpha: \alg{R}_{E/F} \alg{G}_{\add}(F) = E \rightarrow \tilde G, \mbox{ or } \tilde e_\alpha : \alg{R}_{L/F} \alg{J}_{E/L}(F) = \alg{J}_{E/L}(L) \rightarrow \tilde G.$$
Depending on whether $2 \alpha$ is not a root or $2 \alpha$ is a root, we find elements:
$$\tilde e_\alpha(u) \in \tilde G, \mbox{ for } u \in E, \mbox{ or } \tilde e_\alpha(p, \ell), \mbox{ for } (p, \ell) \in E^2, \ell+ \ell^\sigma + p p^\sigma = 0.$$
In the latter case, choose $\theta \in E$ as before, and define $\tilde e_{2 \alpha}(\ell) = \tilde e_\alpha(0, \ell \theta)$.  Then $\tilde e_{2 \alpha}$ coincides with the homomorphism obtained from the unipotent splitting applied to $\alg{e}_{2 \alpha}: \alg{R}_{L/F} \alg{G}_{a, L} \rightarrow \alg{G}$.

Define lifts of $n_\alpha$, for $z \in E^\times$ or for $(c,d) \in \alg{J}_{E/L}^\ast(L)$ by:
\begin{eqnarray*}
\tilde n_\alpha(z) & = & \tilde e_\alpha(z) \tilde e_{-\alpha}(-z^{-1}) \tilde e_\alpha(z), \mbox{ or } \\
\tilde n_\alpha(c,d) & = & \tilde e_{-\alpha}(-c d^{-1}, d^{-\sigma}) \tilde e_\alpha(c,d) \tilde e_{-\alpha}(-c d^{-\sigma}, d^{-\sigma}).
\end{eqnarray*}
In the latter case, define also $\tilde n_{2 \alpha}(\ell) = \tilde n_\alpha(0, \ell \theta)$.  In other terms,
$$\tilde n_{2 \alpha}(\ell) = \tilde e_{-2 \alpha} \left( \frac{\ell^{-1}}{\theta \theta^{\sigma}} \right) \tilde e_{2 \alpha}(\ell) \tilde e_{-2 \alpha} \left( \frac{\ell^{-1}}{ \theta \theta^{\sigma}} \right).$$

Finally, define $\tilde h_\alpha(z)$ for $z \in E^\times$ or $\tilde h_{2 \alpha}(\ell)$ for $\ell \in L^\times$ by
$$\tilde h_\alpha(z) = \tilde n_\alpha(z) \tilde n_\alpha(-1), \mbox{ or } \tilde h_{2 \alpha}(\ell) = \tilde n_{2 \alpha}(\ell) \tilde n_{2 \alpha}(-1).$$

\begin{theorem}
If $2 \alpha$ is not a root, then
\begin{eqnarray*}
\tilde n_\alpha(z) \tilde e_\alpha(u) \tilde n_\alpha(z)^{-1} & = & \tilde e_{-\alpha}(-z^{-2} u), \\
\tilde n_\alpha(z) \tilde e_{-\alpha}(u) \tilde n_\alpha(z)^{-1} & = & \tilde e_\alpha(- z^2 u).
\end{eqnarray*}
If $2 \alpha$ is a root, then
\begin{eqnarray*}
\tilde n_\alpha(c, d) \tilde e_\alpha(p,\ell) \tilde n_\alpha(c,d)^{-1} & = & \tilde e_{-\alpha} \left( -\frac{d^\sigma p}{ d^2}, \frac{\ell}{d d^\sigma} \right), \\
\tilde n_\alpha(c,d) \tilde e_{-\alpha}(p,\ell) \tilde n_\alpha(c,d)^{-1} & = & \tilde e_\alpha \left( \frac{d^2 p^\sigma}{d^\sigma}, \ell d d^\sigma \right).
\end{eqnarray*}
In particular, if $2 \alpha$ is a root, then
\begin{eqnarray*}
\tilde n_\alpha(c, d) \tilde e_{2 \alpha}(\ell) \tilde n_\alpha(c,d)^{-1} & = & \tilde e_{-2 \alpha} \left(\frac{\ell}{d d^\sigma} \right), \\
\tilde n_\alpha(c,d) \tilde e_{-2 \alpha}(\ell) \tilde n_\alpha(c,d)^{-1} & = & \tilde e_{2 \alpha} \left(\ell d d^\sigma \right).
\end{eqnarray*}
\end{theorem}
\proof
The proof is essentially the same as the proof of Theorem \ref{ChevalleySteinberg} in the split case.  If $2 \alpha$ is not a root, we consider (for any $z \in E^\times$) the commutative diagram of groups and homomorphisms over $F$:
$$\xymatrix{
\alg{R}_{E/F} \alg{G}_{\add} \ar[rr]^{\alg{e}_\alpha} \ar[d]_{\alg{\Mult}(-z^{-2})} & & \alg{G} \ar[d]^{\alg{\Int}(n_\alpha(z))} \\
\alg{R}_{E/F} \alg{G}_{\add} \ar[rr]^{\alg{e}_{-\alpha}} & & \alg{G}.
}$$
If $2 \alpha$ is a root, we consider (for any $(c,d) \in \alg{J}_{E/L}^\ast(L)$) the commutative diagram of groups and homomorphisms over $F$:
$$\xymatrix{
\alg{R}_{L/F} \alg{J}_{E/L} \ar[rr]^{\alg{e}_\alpha} \ar[d]_{\alg{f}_{c,d}} & & \alg{G} \ar[d]^{\alg{\Int}(n_\alpha(c,d))} \\
\alg{R}_{L/F} \alg{J}_{E/L} \ar[rr]^{\alg{e}_{-\alpha}} & & \alg{G},
}$$
where $\alg{f}_{c,d}$ is the $L$-automorphism of $\alg{J}_{E/L}$ (or $F$-automorphism of $\alg{R}_{L/F} \alg{J}_{E/L}$) given by:
$$f_{c,d}(p, \ell) = \left( -\frac{d^\sigma p}{ d^2}, \frac{\ell}{d d^\sigma} \right).$$

Since $\tilde G$ is a unipotently split central extension of $G$, the above diagrams lift to give the desired relations in $\tilde G$. 
\qed

From this and previous definitions, we find
\begin{corollary}
If $2 \alpha$ is not a root, then
$$\tilde n_\alpha(z) = \tilde e_\alpha(z) \tilde e_{-\alpha}(-z^{-1}) \tilde e_\alpha(z) = \tilde e_{-\alpha}(-z^{-1}) \tilde e_\alpha(z) \tilde e_{-\alpha}(-z^{-1}).$$
If $2 \alpha$ is a root, then
\begin{eqnarray*}
\tilde n_\alpha(c,d) & = &  \tilde e_{-\alpha}(-c d^{-1}, d^{-\sigma}) \tilde e_\alpha(c,d) \tilde e_{-\alpha}(-c d^{-\sigma}, d^{-\sigma}) \\
& = & \tilde e_\alpha(- c^\sigma d^2 d^{-2 \sigma}, d) \tilde e_{-\alpha}(-c d^\sigma d^{-2}, d^{-\sigma}) \tilde e_\alpha(- c^\sigma d d^{-\sigma}, d).
\end{eqnarray*}
In particular, if $2 \alpha$ is a root, then
\begin{eqnarray*}
\tilde n_{2 \alpha}(\ell) & = & \tilde e_{-2 \alpha} \left( \frac{\ell^{-1}}{\theta \theta^\sigma} \right) \tilde e_{2 \alpha}(\ell) \tilde e_{-2 \alpha} \left( \frac{\ell^{-1}}{\theta \theta^\sigma} \right) \\
& = & \tilde e_{2 \alpha}(\ell) \tilde e_{-2 \alpha} \left( \frac{\ell^{-1}}{\theta \theta^\sigma} \right) \tilde e_{2 \alpha}(\ell).
\end{eqnarray*}
\end{corollary}
\proof
The first case is discussed in Section 11.1 of \cite{B-D}, and the second is essentially contained in \cite{Deo}.  We follow the observation of \cite{B-D}, using the fact that $\tilde n_\alpha(z)$ is invariant under $\Int(n_\alpha(z))$, and $\tilde n_\alpha(c,d)$ is invariant under $\Int(n_\alpha(c,d))$.  

In the first case,
\begin{eqnarray*}
\tilde e_\alpha(z) \tilde e_{-\alpha}(-z^{-1}) \tilde e_\alpha(z) & = & \Int(n_\alpha(z)) \left( \tilde e_\alpha(z) \tilde e_{-\alpha}(-z^{-1}) \tilde e_\alpha(z) \right) \\
&  = & \Int(n_\alpha(z)) \tilde e_\alpha(z) \cdot \Int(n_\alpha(z)) \tilde e_{-\alpha}(-z^{-1}) \cdot \Int(n_\alpha(z))  \tilde e_\alpha(z) \\
& = &  \tilde e_{-\alpha}(-z^{-1}) \tilde e_\alpha(z) \tilde e_{-\alpha}(-z^{-1}).
\end{eqnarray*}
The last step follows from the previous theorem.

In the second case,
\begin{eqnarray*}
& &  \tilde e_{-\alpha}(-c d^{-1}, d^{-\sigma}) \tilde e_\alpha(c,d) \tilde e_{-\alpha}(-c d^{-\sigma}, d^{-\sigma}) \\
& = & \Int(n_\alpha(c,d)) \left(  \tilde e_{-\alpha}(-c d^{-1}, d^{-\sigma}) \tilde e_\alpha(c,d) \tilde e_{-\alpha}(-c d^{-\sigma}, d^{-\sigma}) \right),  \\
& = & \tilde e_\alpha(- c^\sigma d^2 d^{-2 \sigma}, d) \tilde e_{-\alpha}(-c d^\sigma d^{-2}, d^{-\sigma}) \tilde e_\alpha(- c^\sigma d d^{-\sigma}, d).
\end{eqnarray*}
\qed

\begin{corollary}
If $2 \alpha$ is not a root, then $\tilde n_\alpha(z) \cdot \tilde n_\alpha(-z) = 1$.  If $2 \alpha$ is a root, then $\tilde n_{2 \alpha}(\ell) \tilde n_{2 \alpha}(-\ell) = 1$.
\end{corollary}
\proof
In the first case, we compute:
\begin{eqnarray*}
\tilde n_\alpha(z) \tilde n_\alpha(-z) & = & \tilde e_\alpha(z) \tilde e_{-\alpha}(-z^{-1}) \tilde e_\alpha(z) \cdot \tilde e_\alpha(-z) \tilde e_{-\alpha}(-z^{-1}) \tilde e_\alpha(-z) \\
& = & 1,
\end{eqnarray*}
where the last equality follows from the fact that $\alg{e}_\alpha$ is a homomorphism from $\alg{R}_{L/F} \alg{G}_{a,L}$ to $\alg{G}$, and so $\tilde e_\alpha$ is a homomorphism from $L$ to $\tilde G$.

In the second case, we compute:
\begin{eqnarray*}
\tilde n_{2 \alpha}(\ell) \tilde n_{2 \alpha}(-\ell) & = & \tilde e_{-2 \alpha} \left( \frac{\ell^{-1}}{\theta \theta^\sigma} \right) \tilde e_{2 \alpha}(\ell) \tilde e_{-2 \alpha} \left( \frac{\ell^{-1}}{\theta \theta^\sigma} \right) \\
& & \tilde e_{-2 \alpha} \left( \frac{-\ell^{-1}}{\theta \theta^\sigma} \right) \tilde e_{2 \alpha}(-\ell) \tilde e_{-2 \alpha} \left( \frac{-\ell^{-1}}{\theta \theta^\sigma} \right) \\
& = & 1.
\end{eqnarray*}
The last equality follows from the fact that $\alg{e}_{2 \alpha}$ is a homomorphism from $\alg{R}_{L/F} \alg{G}_{a,L}$ to $\alg{G}$, and its lift $\tilde e_{2 \alpha}$ is a homomorphism from $L$ to $\tilde G$.
\qed

\begin{corollary}
If $2 \alpha$ is not a root, then
$$\tilde n_\alpha(z) \tilde n_\alpha(v) \tilde n_\alpha(z)^{-1} = \tilde n_\alpha(z^2 v^{-1}).$$
If $2 \alpha$ is a root, then
$$\tilde n_\alpha(c,d) \tilde n_{2 \alpha}(\ell) \tilde n_\alpha(c,d)^{-1} = \tilde n_{2 \alpha} \left( \frac{d d^\sigma}{\theta \theta^\sigma} \ell^{-1} \right).$$
\end{corollary}
\proof
In the first case, we compute:
\begin{eqnarray*}
\Int(n_\alpha(z)) \tilde n_\alpha(v) & = & \Int(n_\alpha(z)) \left( \tilde e_\alpha(v) \tilde e_{-\alpha}(-v^{-1}) \tilde e_\alpha(v) \right), \\
& = & \tilde e_{-\alpha}(-z^{-2} v) \tilde e_\alpha(z^2 v^{-1}) \tilde e_{-\alpha}(-z^{-2} v), \\
& = & \tilde e_\alpha(z^2 v^{-1}) \tilde e_{-\alpha}(-z^{-2} v) \tilde e_\alpha(z^2 v^{-1}), \\
& = & \tilde n_\alpha(z^2 v^{-1}).
\end{eqnarray*}

In the second case, we compute:
\begin{eqnarray*}
\Int(n_\alpha(c,d)) \tilde n_{2 \alpha}(\ell) & = & \Int(n_\alpha(c,d)) \left(  \tilde e_{-2 \alpha} \left( \frac{\ell^{-1}}{\theta \theta^\sigma} \right) \tilde e_{2 \alpha}(\ell) \tilde e_{-2 \alpha} \left( \frac{\ell^{-1}}{\theta \theta^\sigma} \right) \right), \\
& = &  \tilde e_{2 \alpha} \left( \frac{\ell^{-1} d d^\sigma}{\theta \theta^\sigma} \right) \tilde e_{-2 \alpha} \left( \frac{\ell}{d d^\sigma} \right) \tilde e_{2 \alpha} \left( \frac{\ell^{-1} d d^\sigma}{\theta \theta^\sigma} \right), \\
& = &  \tilde e_{-2 \alpha} \left( \frac{\ell}{d d^\sigma} \right) \tilde e_{2 \alpha} \left( \frac{\ell^{-1} d d^\sigma}{\theta \theta^\sigma} \right) \tilde e_{-2 \alpha} \left( \frac{\ell}{d d^\sigma} \right), \\
& = & \tilde n_{2 \alpha} \left( \frac{d d^\sigma}{\theta \theta^\sigma} \ell^{-1} \right).
\end{eqnarray*}
\qed

Note that the functions $\tilde h_\alpha: E^\times \rightarrow \tilde G$ or $\tilde h_{2 \alpha}: L^\times \rightarrow \tilde G$ (when $\alpha \in \Phi_2$ or $\alpha \not \in \Phi_2$, respectively) are not necessarily homomorphisms.  Rather, as in the split case, there is a 2-cocycle $\sigma_\alpha \in Z^2(E^\times, \mu)$ or $\sigma_{2 \alpha} \in Z^2(L^\times, \mu)$:
\begin{eqnarray*}
\sigma_\alpha(v_1, v_2) & = & \tilde h_\alpha(v_1) \tilde h_\alpha(v_2) \tilde h_\alpha(v_1 v_2)^{-1}, \mbox{ or } \\
\sigma_{2 \alpha}(\ell_1, \ell_2) & = & \tilde h_{2 \alpha}(\ell_1) \tilde h_{2 \alpha}(\ell_2) \tilde h_{2 \alpha}(\ell_1 \ell_2)^{-1}.
\end{eqnarray*}
A simple computation demonstrates that 
\begin{eqnarray*}
\sigma_\alpha(E^\times, 1) = \sigma_\alpha(1, E^\times) = \{ 1 \}, \\
\sigma_{2 \alpha}(L^\times, 1) = \sigma_{2 \alpha}(1, L^\times) = \{ 1 \}.
\end{eqnarray*}

\begin{corollary}
\label{WeylCover}
If $2 \alpha$ is not a root, then
$$\tilde n_\alpha(z) \tilde h_\alpha(v) \tilde n_\alpha(z)^{-1} =  \tilde h_\alpha(v^{-1}) \cdot \sigma_\alpha(v^{-1}, z^2)^{-1}.$$
If $2 \alpha$ is a root, then
$$\tilde n_\alpha(c,d) \tilde h_{2 \alpha}(\ell) \tilde n_\alpha(c,d)^{-1} = \tilde h_{2 \alpha}(\ell^{-1}) \cdot \sigma_{2 \alpha} \left( \ell^{-1}, \frac{d d^\sigma}{\theta \theta^\sigma} \right)^{-1}.$$
\end{corollary}
\proof
In the first case, we compute:
\begin{eqnarray*}
\Int(n_\alpha(z)) \tilde h_\alpha(v) & = & \Int(n_\alpha(z)) \left( \tilde n_\alpha(v) \tilde n_\alpha(-1) \right), \\
& = &  \Int(n_\alpha(z)) \tilde n_\alpha(v) \cdot  \Int(n_\alpha(z))  \tilde n_\alpha(-1), \\
& = & \tilde n_\alpha(z^2 v^{-1}) \tilde n_\alpha(- z^2), \\
& = & \tilde n_\alpha(z^2 v^{-1}) \tilde n_\alpha(-1) \tilde n_\alpha(-1)^{-1} \tilde n_\alpha(z^2)^{-1}, \\
& = & \tilde h_\alpha(z^2 v^{-1}) \tilde h_\alpha(z^2)^{-1}, \\
& = & \tilde h_\alpha(v^{-1}) \cdot \sigma_\alpha(v^{-1}, z^2)^{-1}.
\end{eqnarray*}

In the second case, we compute:
\begin{eqnarray*}
\Int(n_\alpha(c,d)) \tilde h_{2 \alpha}(\ell) & = & \Int(n_\alpha(c,d)) \left( \tilde n_{2 \alpha}(\ell) \tilde n_{2 \alpha}(-1) \right) \\
& = & \Int(n_\alpha(c,d)) \tilde n_{2 \alpha}(\ell) \cdot \Int(n_\alpha(c,d)) \tilde n_{2 \alpha}(-1) \\
& = & \tilde n_{2 \alpha} \left( \frac{d d^\sigma}{\theta \theta^\sigma} \ell^{-1} \right)  \tilde n_{2 \alpha} \left( - \frac{d d^\sigma}{\theta \theta^\sigma} \right), \\
& = & \tilde n_{2 \alpha} \left( \frac{d d^\sigma}{\theta \theta^\sigma} \ell^{-1} \right) \tilde n_{2 \alpha}(-1) \tilde n_{2 \alpha}(-1)^{-1} \tilde n_{2 \alpha} \left(\frac{d d^\sigma}{\theta \theta^\sigma} \right)^{-1}, \\
& = & \tilde h_{2 \alpha} \left( \frac{d d^\sigma}{\theta \theta^\sigma} \ell^{-1} \right) \tilde h_{2 \alpha} \left(\frac{d d^\sigma}{\theta \theta^\sigma} \right)^{-1} \\
& = & \tilde h_{2 \alpha}(\ell^{-1}) \cdot \sigma_{2 \alpha} \left( \ell^{-1}, \frac{d d^\sigma}{\theta \theta^\sigma} \right)^{-1}.
\end{eqnarray*}

\qed

\begin{remark}
Later, we will work in a situation where it is guaranteed that $\sigma_\alpha$ and $\sigma_{2 \alpha}$ are {\em bimultiplicative}, which simplifies the above proposition.  This bimultiplicativity might already follow from relations proven above (cf.~Steinberg \cite{Ste} and Deodhar \cite{Deo}).
\end{remark}

\section{Brylinski-Deligne Extensions}

Let $\alg{G}$ be an affine algebraic group over a field $F$.  In \cite{B-D}, Brylinski and Deligne study a class of central extensions of $\alg{G}$, which are not algebraic groups, but still have algebraic origin.  This entire section can be seen as a review of the results of \cite{B-D}; the few original results in this section are immediate consequences of the deep and beautiful results of Brylinski and Deligne.

\begin{definition}
A central extension of $\alg{G}$ by $\alg{K}_2$ (over $F$) is a central extension of $\alg{G}$ by $\alg{K}_2$, in the category of sheaves of groups on the big Zariski site (of schemes of finite type) over $F$.  Such an extension $\alg{G}'$ is written in a short exact sequence (of sheaves of groups):
$$1 \rightarrow \alg{K}_2 \rightarrow \alg{G}' \rightarrow \alg{G} \rightarrow 1.$$
We write $\Cat{CExt}(\alg{G}, \alg{K}_2)$ for the category of central extensions of $\alg{G}$ by $\alg{K}_2$.
\end{definition}

\begin{remark}
We could be more careful, and refer to a central extension of $\alg{G}$ by $\alg{K}_2$ as a triple $(\alg{G}', \alg{p}, \alg{\iota})$, as in the first section.  But we sacrifice this care in favor of abbreviated notation in this section.
\end{remark}

\begin{remark}
One could also work with central extensions of $\alg{G}$ by $\alg{K}_1$ over $F$, in the category of sheaves of groups on the big Zariski site over $F$.  If one works in the big Zariski site, whose objects are {\em smooth} schemes of finite type over $F$, such central extensions are precisely the central extensions of $\alg{G}$ by $\alg{G}_{\mult}$, studied in the first section.
\end{remark}

A central extension $\alg{G}'$ of $\alg{G}$ by $\alg{K}_2$ yields, for any finitely-generated $F$-algebra $A$, a left-exact sequence of groups:
$$0 \rightarrow \alg{K}_2(A) \rightarrow \alg{G}'(A) \rightarrow \alg{G}(A),$$
in which $\alg{K}_2(A)$ is a subgroup of the center of $\alg{G}'(A)$.  When $A \rightarrow B$ is a morphism of $F$-algebras, there is an obvious commutative diagram, whose rows are left-exact sequences as above; this defines a functor from the category of $F$-algebras to the category of left-exact sequences of groups.

When $L$ is a {\em field} containing $F$, the resulting left-exact sequence is also right exact:
$$0 \rightarrow \alg{K}_2(L) \rightarrow \alg{G}'(L) \rightarrow \alg{G}(L) \rightarrow 0.$$
This arises from the vanishing of $H_{\Zar}^1(L, \alg{K}_2)$ -- a Zariski topology version of Hilbert's Theorem 90 for $\alg{K}_2$.  When $L$ is a Galois extension of $F$, the above short exact sequence is Galois-equivariant.  However, beware that the set of $\Gal(L/F)$-fixed points of $\alg{K}_2(L)$ is often not equal to $\alg{K}_2(F)$.

The category $\Cat{CExt}(\alg{G}, \alg{K}_2)$ is a strictly commutative Picard groupoid, whose structure (at least when $\alg{G}$ is a connected reductive group or parabolic subgroup thereof) is the focus of \cite{B-D}.

\subsection{Central extensions of split unipotent groups by $\alg{K}_2$}

We begin by recalling some of the more basic results of the article \cite{B-D} and their consequences.  The first result provides unipotent splittings:
\begin{proposition}[Prop. 3.2 of \cite{B-D}]
If $\alg{U}$ is a split unipotent group over $F$, then every central extension $1 \rightarrow \alg{K}_2 \rightarrow \alg{U}' \rightarrow \alg{U} \rightarrow 1$ splits uniquely.
\end{proposition}
In other words, the groupoid $\Cat{CExt}(\alg{U}, \alg{K}_2)$ is equivalent to the groupoid with one object and one morphism.

\begin{corollary}
Let $\alg{G}'$ be a central extension of a group scheme $\alg{G}$ by $\alg{K}_2$, over $F$.  Let $G' = \alg{G}'(F)$ and $K_2 = \alg{K}_2(F)$.  Then the previous proposition endows the central extension $G'$ of $G$ by $K_2$ with a unipotent splitting.
\end{corollary}
\proof
Let $\alg{p}: \alg{G}' \rightarrow \alg{G}$ denote the projection homomorphism.  The existence and uniqueness of splitting in the previous proposition yields the following:
\begin{enumerate}
\item
For each $\alg{\eta}: \alg{U} \hookrightarrow \alg{G}$, an embedding of a split unipotent subgroup, the previous proposition yields a unique $\alg{\eta}': \alg{U} \rightarrow \alg{G}'$, satisfying $\alg{p} \circ \alg{\eta}' = \alg{\eta}$:
$$\xymatrix{\alg{U} \ar@/_/[rr]_{\alg{\eta}} \ar[r]^{\alg{\eta}'} & \alg{G}' \ar[r]^{\alg{p}} & \alg{G}.}$$
\item
For every pair $\alg{U}_1$, $\alg{U}_2$ of split unipotent groups, and commutative diagram of closed embeddings:
$$\xymatrix{
\alg{U_1} \ar@/_/[rr]_{\alg{\eta}_1} \ar[r]^{\alg{f}} & \alg{U}_2 \ar[r]^{\alg{\eta}_2} & \alg{G},}$$
the uniqueness in the previous proposition gives a commutative diagram:
$$\xymatrix{
\alg{U}_1 \ar@/_/[rr]_{\alg{\eta}_1'} \ar[r]^{\alg{f}} & \alg{U}_2 \ar[r]^{\alg{\eta}_2'} & \alg{G}', }$$
\item
For each closed embedding of a split unipotent group, $\eta: \alg{U} \hookrightarrow \alg{G}$, and each element $g \in G$, the previous proposition implies that the following diagram commutes:
$$\xymatrix{
\alg{U} \ar[r]^{\alg{\eta}'} \ar[dr]_{ [\alg{\Int}(g) \circ \alg{\eta}]'} & \alg{G}' \ar[d]^{\alg{\Int}(g)} \\
& \alg{G}'
}.$$
\end{enumerate}

Taking $F$-points in each of the three commutative diagrams yields the unipotent splitting.
\qed

\subsection{Central extensions of tori by $\alg{K}_2$}

When $\alg{T}$ is a split torus over $F$, with characters $X = X(\alg{T})$ and cocharacters $Y = Y(\alg{T})$, the category of central extensions of $\alg{T}$ by $\alg{K}_2$ is described in Section 3 of \cite{B-D}.  We describe their result and construction below.
\begin{proposition}[Prop. 3.11 of \cite{B-D}]
Let $\alg{T}$ be a split torus over $F$.  The category of central extensions of $\alg{T}$ by $\alg{K}_2$ is equivalent to the category of pairs $(Q, E)$, where $Q \in \Sym^2 (X)$ is a $\ZZ$-valued quadratic form on the cocharacter lattice $Y$ of $\alg{T}$, and $E$ is a central extension of $Y$ by $F^\times$ (as groups), whose commutator pairing $\Comm: \bigwedge^2 Y \rightarrow F^\times$ is given by:
$$\Comm(y_1, y_2) = (-1)^{B_Q(y_1, y_2)}, \mbox{ where }$$
$$B_Q(y_1, y_2) = Q(y_1 + y_2) - \left( Q(y_1) + Q(y_2) \right)$$
is the symmetric bilinear form associated to $Q$.

When $\alg{T}$ is split over a finite Galois extension $L/F$, the category of central extensiosn of $\alg{T}$ by $\alg{K}_2$ (over $F$) is equivalent to the category of pairs $(Q, E)$ where $Q \in \Sym^2(X)$ as before and $E$ is a $\Gal(L/F)$-equivariant central extension of $Y$ by $L^\times$, satisfying the conditions above.
\end{proposition}
To clarify, consider an extension $E$ of $Y$ by $F^\times$:
$$1 \rightarrow F^\times \rightarrow E \rightarrow Y \rightarrow 1.$$
First, note that $F^\times = \alg{K}_1(F)$, and $\ZZ = \alg{K}_0(F)$ -- this central extension is very much like the extension of $Y$ by $\ZZ$ considered in Theorem \ref{CEByGm}, but with $\alg{G}_{\mult} = \alg{K}_1$ (over a field) replaced by $\alg{K}_2$, and $\alg{K}_0$ replaced by $\alg{K}_1$.  Second, note that such a central extension yields a commutator pairing by defining:
$$\Comm(y_1, y_2) = e_1 e_2 e_1^{-1} e_2^{-1},$$
for any $e_1, e_2 \in E$ projecting to $y_1$ and $y_2$ respectively.

The previous proposition implies that for a split torus $\alg{T}$, $\Cat{CExt}(\alg{T},\alg{K}_2$) is a groupoid whose isomorphism classes are parameterized by $\Sym^2(X)$, and all of whose automorphism groups are isomorphic to $\Hom(Y, F^\times) = \alg{T}^\vee(F)$, where $\alg{T}^\vee = Spec(F[Y])$ is the dual torus to $\alg{T}$.

Since it would otherwise be completely mysterious, we describe the functor from $\Cat{CExt}(\alg{T}, \alg{K}_2)$ to $\Cat{CExt}(Y, F^\times)$ (the latter being central extensions in the category of groups) explicitly, in three steps.  Begin with a central extension $\alg{T}'$ of $\alg{T}$ by $\alg{K}_2$.
\begin{enumerate}
\item
Taking points over the Laurent series field $F((\tau))$, one gets a central extension of groups:
$$1 \rightarrow \alg{K}_2(F((\tau))) \rightarrow \alg{T}'(F((\tau))) \rightarrow \alg{T}(F((\tau))) \rightarrow 1.$$
\item
Sending a cocharacter $y \in Y$ to $y(\tau) \in \alg{T}(F((\tau)))$ gives us an embedding of groups $Y \hookrightarrow \alg{T}(F((\tau)))$, and allows us to pull back this central extension:
$$1 \rightarrow \alg{K}_2(F((\tau))) \rightarrow Y' \rightarrow Y \rightarrow 1.$$
\item
Pushing forward via the tame symbol (see Definition \ref{TameSymbol}) $\tame_{F((\tau))}:  \alg{K}_2(F((\tau))) \rightarrow \alg{K}_1(F) = F^\times$ yields a central extension
$$1 \rightarrow F^\times \rightarrow E \rightarrow Y \rightarrow 1.$$
\end{enumerate}
This construction is described in Section 3.10 and Remark 3.12 of \cite{B-D}.  It generalizes naturally to the nonsplit case.

All central extensions of a split torus $\alg{T}$ by $\alg{K}_2$ are ``incarnated'' by a 2-cocycle of algebraic origin.  The following construction is described in Sections 3.9-10 of \cite{B-D}:  begin with a quadratic form $Q \in \Sym^2(X)$.  Choose any representative bilinear form 
$$D = \sum_i (x_1^i \otimes x_2^i)  \in X \otimes_\ZZ X,$$ 
projecting to $Q$, i.e., $Q(y) = D(y,y)$ for all $y \in Y$.  Then, one may define a central extension $\alg{T}'$ of $\alg{T}$ by $\alg{K}_2$, endowed with a trivialization of the $\alg{K}_2$-torsor $\alg{T}'$ over $\alg{T}$, whose 2-cocycle is given by a finite product
$$\sigma(t_1, t_2) = \prod_i \{ x_1^i(t_1), x_2^i (t_2) \}.$$
In other words, if $L$ is a field containing $F$, then we may identify $\alg{T}'(L)$ as a set with $\alg{T}(L) \times \alg{K}_2(L)$; the group law on $\alg{T}'(L)$ is given by the usual group law on the central subgroup $\alg{K}_2(L)$ (written multiplicatively), and the following ``twisted'' multiplication:
$$(t_1,1) \cdot (t_2,1) = (t_1 t_2, \sigma(t_1, t_2)) = \left( t_1 t_2, \prod_i \{ x_1^i(t_1), x_2^i(t_2) \}_L \right).$$
Here, we write $\{ \cdot, \cdot \}_L$ for the symbol from $L^\times \times L^\times$ to $\alg{K}_2(L)$.  This construction yields a central extension whose isomorphism class has parameter $Q$.

Conversely, if one is given a central extension $\alg{T}'$ of a split torus $\alg{T}$ by $\alg{K}_2$, then there exists a section $\alg{j}: \alg{T} \rightarrow \alg{T}'$ of the underlying $\alg{K}_2$-torsor, which satisfies $\alg{j}(1) = 1$.  Such a section gives an algebraic cocycle $\alg{\sigma}: \alg{T} \times \alg{T} \rightarrow \alg{K}_2$ (a section of the Zariski sheaf $\alg{K}_2$ over $\alg{T} \times \alg{T}$), given at the level of points by
$$\sigma(t_1, t_2) = j(t_1) j(t_2) j(t_1 t_2)^{-1}.$$
Since $\alg{\sigma}$ is trivial on $\{ 1 \} \times \alg{T}$ and $\alg{T} \times \{ 1 \}$, Corollary 3.7 of \cite{B-D} ensures that the section $\alg{\sigma}$ is bimultiplicative.  This result is the $\alg{K}_2$-analogue of the $\alg{G}_{\mult}$-result proven at the end of Theorem \ref{CEByGm} and discussed in the subsequent remark.  Any such bimultiplicative section has the form
$$\sigma(t_1, t_2) = \prod_i \{ x_1^i(t_1), x_2^i (t_2) \},$$
for some element $\sum_i (x_1^i \otimes x_2^i) \in X \otimes_\ZZ X$ as above.  While this element of $X \otimes_\ZZ X$ is not uniquely determined by $\alg{\sigma}$, it is uniquely determined up to the subgroup $\bigwedge^2 X \subset X \otimes_\ZZ X$; this yields a well-defined element $Q \in \Sym^2 X$, from any central extension $\alg{T}'$ of $\alg{T}$.    

It will be important to understand the central extensions of tori by $\alg{K}_2$, in one nonsplit situation.  
\begin{example}
\label{RSExamp}
Consider a separable extension $L/F$, and the Weil restriction of scalars $\alg{T} = \alg{R}_{L/F} \alg{G}_{\mult}$.  Let $E$ be a normal closure of $L$, $\Gamma = \Gal(E/F)$, and $I$ the set of $F$-algebra embeddings from $L$ into $E$.  Thus $I$ is identified as a $\Gamma$-set with $\Gamma / \Gal(E/L)$.  The characters and cocharacters of $\alg{T}$ can be identified as $X = \ZZ^I$, $Y = \ZZ^I$.

Any central extension of $\alg{T}$ by $\alg{K}_2$ is obtained via descent from a central extension of $\alg{G}_{\mult}^I$ by $\alg{K}_2$ over $E$; its isomorphism class depends on a quadratic form $Q: Y = \ZZ^I \rightarrow \ZZ$, which is $\Gamma$-invariant.  Such a quadratic form can be represented by a $\Gamma$-invariant bilinear form $D: Y \otimes Y \rightarrow \ZZ$, in the sense that
$$Q(y) = D(y,y), \quad D(y_1, y_2) = D(\gamma y_1, \gamma y_2),$$
for all $y, y_1, y_2 \in Y$ and $\gamma \in \Gamma$.  This can be seen by direct construction, or by a cohomological argument since $\bigwedge^2 Y$ is a sum of induced modules over $\ZZ[\Gamma]$.  

Now consider the canonically embedded $\alg{G}_{\mult} \hookrightarrow \alg{T}$ over $F$.  It corresponds to the diagonal embedding $\Delta: \ZZ \hookrightarrow Y = \ZZ^I$.  We find that
$$D(\Delta(m), \Delta(n)) = (\# I) \cdot D( (m,0, \ldots, 0), \Delta(n)) \in (\# I) \cdot \ZZ.$$

If $\alg{G}_{\mult}'$ is the central extension of $\alg{G}_{\mult}$ by $\alg{K}_2$, obtained by pulling back $\alg{T}'$ via $\Delta$, then we find that $\alg{G}_{\mult}'$ is incarnated by a cocycle of the form:
$$\sigma(\Delta(z_1), \Delta(z_2)) = \{ z_1, z_2 \}^{\# I \cdot d},$$
where $d = D((1,0,\ldots, 0), (1,\ldots, 1)) \in \ZZ$.

The central extension $\alg{T}'$ is incarnated by some $\Gamma$-invariant bimultiplicative cocycle $\alg{\sigma}: \alg{T} \times \alg{T} \rightarrow \alg{K}_2$ (This follows from Theorem 2.1 of \cite{B-D}), defined over $F$, extending the above cocycle on $\Delta(\alg{G}_{\mult})$.  Furthermore, we can compute:
$$\sigma(t_1, \Delta(z_2) )^{\# I} = \sigma( N_{L/F}(t_1), \Delta(z_2) ) = \{ N_{L/F} t_1, z_2 \}^{\# I \cdot d }.$$
Since the group of bimultiplicative sections of $\alg{K}_2$ over $\alg{T} \times \alg{G}_{\mult}$ is torsion-free (via Corollary 3.7 of \cite{B-D} and descent), this implies that
$$\sigma(t_1, \Delta(z_2)) = \{ N_{L/F} t_1, z_2 \}^d.$$
\end{example}

\subsection{Central extensions of Chevalley groups by $\alg{K}_2$}

Let $\alg{G}$ be a split semisimple simply-connected group over $F$, with $F$-split maximal torus $\alg{T}$.  From Brylinski and Deligne, we recall the following
\begin{theorem}[Special case 7.3(i) of \cite{B-D}]
The isomorphism classes in the category of central extensions of $\alg{G}$ by $\alg{K}_2$ are in bijection with the $W$-invariant $\ZZ$-valued quadratic forms $Q: Y \rightarrow \ZZ$, i.e. the elements $Q \in \Sym^2(X)^W$.  There are no nontrivial automorphisms in the category of central extensions of $\alg{G}$ by $\alg{K}_2$.  
\end{theorem}
\begin{remark}
This parameterization is compatible (Compatibility 4.9 of \cite{B-D}) with the parameterization for split tori -- when $\alg{T}$ is a split torus in $\alg{G}$, and $\alg{G}'$ is a central extension of $\alg{G}$ by $\alg{K}_2$, one may pull back $\alg{G}'$ to get a central extension $\alg{T}'$ of $\alg{T}$ by $\alg{K}_2$.  The invariant $Q \in \Sym^2(X)^W$ associated to $\alg{G}'$ is equal to the invariant $Q \in \Sym^2(X)$ associated in the previous section to $\alg{T}'$.
\end{remark}

For comparison, recall that when $\alg{G}$ is a simply-connected split semisimple group over the perfect field $\FF$, there is exactly one isomorphism class of central extensions of $\alg{G}$ by $\alg{G}_{\mult}$ -- the class of the split extension $\alg{G} \times \alg{G}_{\mult}$ -- and the split extension has no nontrivial automorphisms.  Brylinski and Deligne prove that central extensions of $\alg{G}$ by $\alg{K}_2$ are slightly more complicated, in that there are numerous isomorphism classes, but the category is still ``rigid'' -- objects have no non-identity automorphisms.

When $\alg{G}$ is simply-connected, split, and {\em almost simple}, the set $\Sym^2(X)^W$ can be identified with $\ZZ$.  This follows from an observation about root systems:
\begin{proposition}
\label{UniqueQuadratic}
For a simple root datum $(X, \Phi, Y, \Phi^\vee)$, there is a unique integer-valued Weyl-invariant quadratic form $Q_1: Y \rightarrow \ZZ$ satisfying the identity $Q_1(\alpha^\vee) = 1$ for all short coroots (those coroots associated to long roots) $\alpha^\vee \in \Phi^\vee$.  Moreover $Sym^2(X)^W = \ZZ \cdot Q_1$.
\end{proposition}
\proof
There is a unique, up to scaling, $\QQ$-valued quadratic invariant polynomial for the reflection representation of a finite irreducible Coxeter group.  Hence if $Q, Q'$ are two integer-valued, Weyl-invariant quadratic forms on $Y$, then $Q' = q \cdot Q$ for some rational number $q$.  

Checking case-by-case (see the following example), we find that if $Q$ is a $\QQ$-valued Weyl-invariant quadratic form on $Y$, and $Q(\alpha^\vee) = 1$ for any short coroot $\alpha^\vee$, then $Q$ is $\ZZ$-valued.  The result follows immediately.
\qed

\begin{example}
In a simply-laced simple root system, every coroot is short, and $Q_1$ takes the value $1$ on every coroot.  In types $\Type{B}$, $\Type{C}$, and $\Type{F}_4$, $Q_1$ takes the value $1$ on short coroots and $2$ on long coroots.  In type $\Type{G}_2$, $Q_1$ takes the value $1$ on short coroots and $3$ on long coroots.
\end{example}

By Proposition 4.15 of \cite{B-D}, when $\alg{G}$ is an almost-simple simply-connected split Chevalley group, the central extension of $\alg{G}$ by $\alg{K}_2$ corresponding to $Q_1$ coincides (upon taking $F$-points) with the universal central extension (outside of type $\Type{C}$) studied by Steinberg \cite{Ste} and Matsumoto \cite{Mats}.

\begin{Notation}
\label{CExtFromQ}
Let $\alg{G}$ be a simply-connected split semisimple group over $F$, and $Q \in \Sym^2(X)^W$ a quadratic form on $Y$.  Let $\alg{T}$ be a $F$-split maximal torus in $\alg{G}$.  Let $\alg{G}_Q'$ be the associated (unique up to unique isomorphism) central extension of $\alg{G}$ by $\alg{K}_2$.  Let $\alg{T}_Q'$ be the resulting central extension of $\alg{T}$ by $\alg{K}_2$.  Define $E_Q$ to be the resulting central extension of $Y$ by $F^\times$ (from the previous section).
\end{Notation}

In other words, the data of a Weyl-invariant quadratic form on $Y$ -- for a split simply-connected semisimple group -- yields a central extension of $Y$ by $F^\times$.  This object of $\Cat{CExt}(Y, F^\times)$ is characterized up to unique isomorphism, in Section 11 of \cite{B-D}. 

\subsection{Simply-connected semisimple groups}

Let $\alg{G}$ be a simply-connected semisimple group over $F$ (no longer necessarily split).  Let $\alg{T}$ be a maximal torus in $\alg{G}$, defined over $F$, and $Y = Y(\alg{T})$ its cocharacter lattice.  Let $F^{\sep}$ be a separable closure of $F$, and $\Gamma = \Gal(F^{\sep} / F)$, so that $Y$ is naturally a $\ZZ[\Gamma]$-module.  Let $\alg{W} = \alg{W}(\alg{G}, \alg{T})$ denote the Weyl group, a finite \'etale group over $F$.  Let $W$ be the geometric points of $\alg{W}$, viewed as a group with action of $\Gamma$.    

We mention a number of results here, which are straightforward consequences of the main result of \cite{B-D}.
\begin{theorem}[Theorem 7.2 and Special case 7.3(i) of \cite{B-D}]
Let $\alg{G}$ be simply-connected and semisimple.  Then $\Cat{CExt}(\alg{G}, \alg{K}_2)$ is a rigid groupoid -- between any two objects there is at most one morphism.  
The isomorphism classes of $\Cat{CExt}(\alg{G}, \alg{K}_2)$ are in natural bijection with the set of $\Gamma$ and $W$ invariant quadratic forms on $Y$, i.e., the central extensions of $\alg{G}$ by $\alg{K}_2$ are classified up to unique isomorphism by elements of $(\Sym^2 X)^{\Gamma \ltimes W}$.
\end{theorem}

\begin{corollary}
Let $\alg{G} = \alg{G}_1 \times \alg{G}_2$ be a product of simply-connected and semisimple groups over $F$.  Then, there is a natural equivalence of rigid groupoids 
$$\Cat{CExt}(\alg{G}, \alg{K}_2) \isom \Cat{CExt}(\alg{G}_1, \alg{K}_2) \times \Cat{CExt}(\alg{G}_2, \alg{K}_2).$$
\end{corollary}
\proof
Choosing $F$-tori $\alg{T}_1$ and $\alg{T}_2$ in $\alg{G}_1$ and $\alg{G}_2$ respectively, with resulting Weyl groups $W_1, W_2$, the central extensions of $\alg{G}$ by $\alg{K}_2$ are classified up to unique isomorphism by elements of $(\Sym^2 (X_1 \oplus X_2))^{\Gamma \ltimes (W_1 \times W_2)}$.  Such a quadratic form decomposes as a sum of elements of $\Sym^2(X_1)^{\Gamma \ltimes W_1}$, $\Sym^2(X_2)^{\Gamma \ltimes W_2}$, and $(X_1 \otimes X_2)^{\Gamma \ltimes (W_1 \times W_2)}$.  The latter ``cross-terms'' must vanish, since for every coroot $\alpha_1^\vee \in Y_1$ (and $Y_1$ is generated by coroots for $\alg{G}_1$ with respect to $\alg{T}_1$), the reflection $w = s_{\alpha_1^\vee}$ satisfies $w \alpha_1^\vee = - \alpha_1^\vee$ and $w \alpha_2^\vee = \alpha_2^\vee$ for all $\alpha_2^\vee \in \Phi_2^\vee$.  No nonzero element of $X_1 \otimes X_2$ (viewed as a bilinear form on $Y_1 \times Y_2$) can be invariant under all such reflections.  

Hence we find the decomposition 
$$(\Sym^2 (X_1 \oplus X_2))^{\Gamma \ltimes (W_1 \times W_2)} = \Sym^2(X_1)^{\Gamma \ltimes W_1} \oplus \Sym^2(X_2)^{\Gamma \ltimes W_2}.$$
\qed

The classification of central extensions of absolutely almost simple groups by $\alg{K}_2$ follows from Proposition \ref{UniqueQuadratic}:
\begin{proposition}
Suppose that $\alg{G}$ is absolutely almost simple, and simply-connected, semisimple as before.  Then the central extensions of $\alg{G}$ by $\alg{K}_2$ are classified up to unique isomorphism by elements of $\ZZ$.  Namely, there is a unique $W$-invariant quadratic form $Q$ on $Y$, such that $Q(\alpha^\vee) = 1$ for every short coroot $\alpha^\vee \in Y$, and every integer multiple of this quadratic form is $\Gamma$-invariant.
\end{proposition}
\proof
In the split case, we have already mentioned the uniqueness of such a quadratic form.  Since the action of $\Gamma$ on $Y$ must send short roots to short roots, it follows that this quadratic form and its integer multiples are $\Gamma$-invariant.
\qed

The almost simple over $F$ case follows from the previous two:
\begin{corollary}
Suppose that $\alg{G}$ is simply-connected, semisimple, and almost simple over $F$ (not necessarily absolutely almost simple).  Then the central extensions of $\alg{G}$ by $\alg{K}_2$ are classified up to unique isomorphism by elements of $\ZZ$.  
\end{corollary}
\proof
Such a group $\alg{G}$ is isomorphic to $\alg{R}_{L/F} \alg{H}$, for some absolutely simple, simply-connected, semisimple group $\alg{H}$ over $L$, with $L$ a finite separable extension of $F$.  An object of $\Cat{CExt}(\alg{G}, \alg{K}_2)$ is determined by an object of $\Cat{CExt}_L(\alg{H}^I, \alg{K}_2)$, endowed with descent data, where $I$ is the set of embeddings of $L$ into a fixed separable closure $F^{\sep}$ of $F$.

An object of $\Cat{CExt}_L(\alg{H}^I, \alg{K}_2)$ is determined by an indexed family of integers $(q_i)_{i \in I}$, by the previous proposition.  For there to exist descent data down to $F$, these integers must be equal; the rigidity of the groupoid $\Cat{CExt}_L(\alg{H}, \alg{K}_2)$ implies the existence (when these integers are equal) and uniqueness of descent data.
\qed

\subsection{Reductive groups}
The main theorem of \cite{B-D} describes, completely and practically, the category $\Cat{CExt}(\alg{G}, \alg{K}_2)$, in the same way that our (easier) Theorem \ref{CEByGm} describes $\Cat{CExt}(\alg{G}, \alg{G}_{\mult})$.  Here is the main theorem of Brylinski and Deligne, given with a bit more detail filled in for the reader:  
\begin{theorem}[Theorem 7.2 and the Introduction of \cite{B-D}]
Let $\alg{G}$ be a connected reductive group over a field $F$.  Let $\alg{T}$ be a maximal torus of $\alg{G}$ defined over $F$, with character group $X$ and cocharacter group $Y$.  Let $L$ be a (finite) Galois extension of $F$ which splits $\alg{T}$.  The category $\Cat{CExt}(\alg{G}, \alg{K}_2)$ is equivalent to the following category of quintuples $(Q,E,p,\iota,\phi)$: the first entry is a quadratic form $Q \in \Sym^2(X)^{\Gamma \ltimes W}$, and $(E, p, \iota)$ is a $\Gamma$-equivariant central extension of $Y$ by $L^\times$:
$$\xymatrix{
1 \ar[r] & L^\times \ar[r]^{\iota} & E \ar[r]^p & Y \ar[r] & 1,}$$
whose commutator pairing $\Comm: \bigwedge^2 Y \rightarrow L^\times$ is required to satisfy $\Comm(y_1, y_2) = (-1)^{B_Q(y_1, y_2)}$.

Let $\alg{f}: \alg{G}_{\sconn} \rightarrow \alg{G}$ be the composition of $\alg{G}_{\sconn} \rightarrow \alg{G}_{der} \rightarrow \alg{G}$, i.e., the simply-connected cover of the derived subgroup of $\alg{G}$, $\alg{T}_{\sconn} = \alg{f}^{-1}(\alg{T})$, and $Y_{\sconn}$ the cocharacter group of $\alg{T}_{\sconn}$.  Let $E_Q$ be the ($\Gamma$-equivariant) central extension of $Y_{\sconn}$ by $L^\times$, associated to the quadratic form $Q$ restricted to $Y_{\sconn}$ by Definition \ref{CExtFromQ}.  The last part of the quintuple, $\phi$, is a $\Gamma$-equivariant morphism from $E_Q$ to $E$ making the following diagram commute:
$$\xymatrix{
0 \ar[r] & L^\times \ar[r] \ar[d]^= & E_Q \ar[r] \ar[d]^{\phi} & Y_{\sconn} \ar[r]  \ar[d]^{f_\ast} & 0 \\
0 \ar[r] & L^\times \ar[r]^{\iota} & E \ar[r]^{p} & Y \ar[r] & 0.
}$$
Morphisms of from a quintuple $(Q_1,E_1, p_1, \iota_1, \phi_1)$ to a quintuple $(Q_2, E_2, p_2, \iota_2, \phi_2)$ exist only when $Q_1 = Q_2$, and in this case are $\Gamma$-equivariant homomorphisms from $E_1$ to $E_2$, making the large but obvious $\Gamma$-equivariant diagram of groups commute. 
\end{theorem}

The above classification is compatible with passage to (standard) Levi subgroups.  Namely, if $\alg{L}$ is a Levi factor of an $F$-parabolic subgroup $\alg{P} \subset \alg{G}$, and $\alg{L}$ contains $\alg{T}$, then one may restrict an object $\alg{G}'$ of $\Cat{CExt}(\alg{G}, \alg{K}_2)$ to obtain an object $\alg{L}'$ of $\Cat{CExt}(\alg{L}, \alg{K}_2)$.  Since $\alg{T}$ is a maximal torus in $\alg{L}$, the data $(Q, E, p, \iota)$ associated to $\alg{G}'$ is also the data associated to $\alg{L}'$.  The last part $\phi$ of the data is slightly more difficult to describe, and is different for $\alg{G}$ and $\alg{L}$.

\section{Tame covers}
In this section, we finally specialize to the case when $F$ is a complete, discretely-valued field with valuation ring $\OO$ and perfect residue field $\FF$.  When convenient, we use a uniformizing element $\varpi$ for $F$, and we normalize the valuation so that $\val(F^\times) = \ZZ$ and $\val(\varpi) = 1$.  When $u \in \OO$, we write $\bar u$ for its reduction in $\FF$.
\begin{definition}
\label{TameSymbol}
The {\em tame symbol} is the homomorphism $\tame_F: \alg{K}_2(F) \rightarrow \FF^\times$ given by:
$$\tame_F( \{x,y \} ) = (-1)^{\val(x) \cdot \val(y)} \cdot \overline{ \left( \frac{y^{\val(x)}}{x^{\val(y)}} \right) },$$
We also write:
$$\{ x,y \}_{\tame} = \tame_F( \{x,y \} ),$$
when the field $F$ is clear from context.
\end{definition}

\begin{definition}
Let $\alg{G}'$ be a central extension of $\alg{G}$ by $\alg{K}_2$.  Consider the central extension $\tilde G$ of $G$ by $\FF^\times$, given as the pushforward in the following diagram:
$$\xymatrix{
1 \ar[r] & \alg{K}_2(F) \ar[d]^{\tame_F} \ar[r] & \alg{G}'(F) \ar[r] \ar[d] & \alg{G}(F) \ar[r] \ar[d]^{=} &  1 \\
1 \ar[r] & \FF^\times \ar[r] & \tilde G \ar[r] & G \ar[r] & 1.
}.$$
In this setting, we say that $\tilde G$ is the {\em tame extension} of $G$ by $\FF^\times$ arising from $\alg{G}'$.
\end{definition}

When $F$ is a nonarchimedean local field, $\tilde G$ is naturally a locally compact group, whose quotient by the finite central subgroup $\FF^\times$ is topologically isomorphic to the locally compact group $G$ (See Construction 10.3 of \cite{B-D}).

\subsection{Tame behavior}

The following properties of the tame symbol are quite useful for computations.  While they can be found in many texts on K-theory, we find Chapter 7.1 of \cite{G-S} (and other chapters) an outstanding reference on the subject.  All facts we use about the tame symbol can be found there.

\begin{proposition}
For all $x,y \in \OO^\times$, $\{x,y \}_{\tame} = 1$.  If $\varpi$ is a uniformizing element of $F$, and $x \in \OO^\times$, then $\{ \varpi, x \}_{\tame} = \bar x$, where $\bar x$ denotes the reduction of $x$ in $\FF^\times$.  Also, $\{ \varpi, \varpi \}_{\tame} = \{ \varpi, -1 \}_{\tame} = \overline{-1}$.
\end{proposition}
\proof
All three claims follow from the definition of the tame symbol.
\qed

\begin{remark}
Often, one works with ``metaplectic groups'' which are obtained (in odd residue characteristic) by pushing the tame cover $\tilde G$ forward further via the local Legendre symbol $\Leg_2: \FF^\times \rightarrow \mu_2$.  Since $\Leg_2(-1)$ (the pushforward of $\{ \varpi, \varpi \}_{\tame}$) depends on the congruence class of $q = \# \FF$ modulo $4$, one finds this sort of arithmetic naturally when working with metaplectic groups.  Quadratic reciprocity arises in the global theory of the metaplectic group.
\end{remark}

The following describes the behavior of the tame symbol when passing to a finite separable field extension.
\begin{proposition}[Chapter 7 of \cite{G-S}]
Let $L$ be a finite separable extension of $F$, of ramification index $e$, with residue field $\LL$.  Then the following diagram commutes, where the map from $\alg{K}_2(F)$ to $\alg{K}_2(L)$ is the homomorphism functorially associated to the inclusion $F \hookrightarrow L$:
$$\xymatrix{
K_2(F) \ar[rr]^{\tame_F} \ar[d] & & \FF^\times \ar[d]^{z \mapsto z^e} \\
K_2(L) \ar[rr]^{\tame_L} & & \LL^\times
}.$$
\end{proposition}
A consequence for tame extensions is the following observation of Brylinski and Deligne:
\begin{corollary}[Proposition 12.9 of \cite{B-D}]
Let $\alg{G}'$ be a central extension of $\alg{G}$ by $\alg{K}_2$.  Let $L$ be an unramified Galois extension of $F$.  Then there is a natural commutative diagram:
$$\xymatrix{
1 \ar[r] & \FF^\times \ar[d] \ar[r] & \tilde G_F \ar[d] \ar[r] & G_F \ar[r]  \ar[d] & 1 \\
1 \ar[r] & \LL^\times \ar[r] & \tilde G_L \ar[r] & G_L \ar[r] & 1,
}$$
in which $G_F = \alg{G}(F)$, $G_L = \alg{G}(L)$, and $\tilde G_F$ and $\tilde G_L$ are the tame extensions arising from pushing forward $\alg{G}'(F)$ and $\alg{G}'(L)$ via $\tame_F$ and $\tame_L$.  Moreover, the top row in this diagram is precisely equal to the $\Gal(L/F)$-fixed points of the bottom row.
\end{corollary}

Finally, we mention the following crucial result of Brylinski and Deligne, on which we elaborate later in this section.
\begin{theorem}[Construction 12.11 of \cite{B-D}]
Let $\alg{\underline G}$ be a smooth group scheme over $\OO$, with generic fibre $\alg{G}$ and special fibre $\alg{\bar G}$.  Let $\alg{G}'$ be a central extension of $\alg{G}$ by $\alg{K}_2$.  Let $\tilde G$ be the resulting tame extension of $G$ by $\FF^\times$.  Let $G^\circ = \alg{\underline{G}}(\OO)$ and let $\tilde G^\circ$ denote the preimage of $G^\circ$ in $\tilde G$.  

Then there exists a central extension $\alg{\bar G}'$ of $\alg{\bar G}$ by $\alg{\bar G}_{\mult}$ over $\FF$, and a commutative diagram with exact rows:
$$\xymatrix{
1 \ar[r] &  \FF^\times \ar[r] \ar[d]^{=} & \tilde G^\circ \ar[r] \ar[d] & G^\circ \ar[d] \ar[r] & 1 \\
1 \ar[r] & \alg{\bar G}_{\mult}(\FF) \ar[r] & \alg{\bar G}'(\FF) \ar[r] & \alg{\bar G}(\FF) \ar[r] & 1
}$$
such that the $\tilde G^\circ$ is obtained via pullback from the central extension in the bottom row. 
\end{theorem}

Now that we have recalled the essential results of \cite{B-D}, we will approach Question 12.13 of \cite{B-D} and describe the central extensions $\alg{\bar G}'$ that arise when $\alg{\underline G}$ comes from Bruhat-Tits theory.  Such a description is crucial for the understanding of the depth-zero representations of $\tilde G$, discussed in work of T. Howard and the author \cite{We2}.
 
\subsection{Tame extensions of split tori}
Let $\alg{S}$ be a split torus over $F$.  Let $\alg{S}'$ be a central extension of $\alg{S}$ by $\alg{K}_2$.  Let $\tilde S$ be the resulting tame extension of $S$ by $\FF^\times$.  Such groups $\tilde S$ and their representation theory were studied in the author's earlier paper \cite{We1}.  Let $X$ be the character group, and $Y$ the cocharacter group of $\alg{S}$.  As a split torus, we may identify $\alg{S}$ with $Spec(F[X])$.  We write $\alg{\underline{S}} = Spec(\OO[X])$ for the canonical model of $\alg{S}$ over $\OO$.  Let $S^\circ = \alg{\underline{S}}(\OO)$ be the resulting subgroup of $S$, and $\tilde S^\circ$ its preimage in the tame extension:
$$\xymatrix{
1 \ar[r] & \FF^\times \ar[r] \ar[d]^{=} & \tilde S^\circ \ar[r] \ar[d] & S^\circ \ar[r] \ar[d] & 1 \\
1 \ar[r] & \FF^\times \ar[r] & \tilde S \ar[r] & S \ar[r] & 1.}$$

Following Corollary 3.7 of \cite{B-D}, there exists a section $\alg{j}$ of the $\alg{K}_2$-torsor $\alg{S}' \rightarrow \alg{S}$ (pointed at the identity).    By Lemma 12.12 of \cite{B-D}, the $\alg{K}_2$-torsor $\alg{S}' \rightarrow \alg{S}$ yields, via the residue map in K-theory, a $\alg{G}_{\mult}$-torsor $\alg{\bar S}' \rightarrow \alg{S}$; this construction yields a functor:
$$Res:  \Cat{CExt}(\alg{S}, \alg{K}_2) \rightarrow \Cat{CExt}(\alg{\bar S}, \alg{\bar G}_{\mult}).$$
Of course, every central extension of $\alg{\bar S}$ by $\alg{\bar G}_{\mult}$ is again a torus -- central extensions are abelian extensions in this situation, but we maintain the notation of the first section here.  
\begin{remark}
It seems likely that this functor can also be seen through the classifications of central extensions (Theorem \ref{CEByGm} and the Main Theorem of \cite{B-D}), sending a central extension $F^\times \rightarrow E \rightarrow Y$ to an extension $\ZZ \rightarrow Y' \rightarrow Y$, by pushing forward via $\val: F^\times \rightarrow \ZZ$.  We have not checked this, however.
\end{remark} 

This construction yields, for any section $\alg{j}$ of the pointed  $\alg{K}_2$-torsor, a section $\alg{\bar j}$ of the corresponding $\alg{G}_{\mult}$-torsor.  If $\alg{\sigma}: \alg{S} \times \alg{S} \rightarrow \alg{K}_2$ is the 2-cocycle associated to $\alg{j}$:
$$\sigma(s_1, s_2) = j(s_1) j(s_2) j(s_1 s_2)^{-1},$$
then $\alg{\sigma}$ is bimultiplicative, and subject to the classification of Corollary 3.7 of \cite{B-D}.  Moreover $\alg{\bar \sigma}$ is the 2-cocycle associated to $\alg{\bar j}$, and is given for all $s_1, s_2 \in S^\circ$ by:
$$\bar \sigma(\bar s_1, \bar s_2) = \bar j(s_1) \bar j(s_2) \bar j(s_1 s_2)^{-1} = \tame(\sigma(s_1, s_2)).$$
Since the tame symbol is trivial on $\OO^\times \times \OO^\times$, we find that $\bar \sigma$ is trivial; therefore $\alg{\bar j}$ is not just a section, but is a {\em splitting} of the cover $\alg{\bar S}' \rightarrow \alg{\bar S}$.  We find that
\begin{proposition}
Every section $\alg{j}$ of the pointed $\alg{K}_2$-torsor $\alg{S}'$ over $\alg{S}$ yields a splitting $\alg{\bar j}$ of the extension
$$1 \rightarrow \alg{\bar G}_{\mult} \rightarrow \alg{\bar S}' \rightarrow \alg{\bar S} \rightarrow 1.$$
By pullback it yields a splitting $j^\circ$ of the tame extension,
$$1 \rightarrow \FF^\times \rightarrow \tilde S^\circ \rightarrow S^\circ \rightarrow 1.$$  
The set of splittings $\alg{\bar j}$, as $\alg{j}$ is allowed to vary, forms a torsor for $X = \Hom(\alg{S}, \alg{G}_{\mult}) = \Hom(\alg{\bar S}, \alg{\bar G}_{\mult})$.
\end{proposition}
\proof
It only remains to check the final claim.  The set of all splittings of $\alg{\bar G}_{\mult} \rightarrow \alg{\bar S}' \rightarrow \alg{\bar S}$ forms a torsor for $X$, so it must be checked that all such splittings arise as reductions from sections $\alg{j}$ of the $\alg{K}_2$ torsor.  For this, observe that for all $x \in X$, and all sections $\alg{j}$ of $\alg{S}' \rightarrow \alg{S}$, one may construct a new section ${}^x \alg{j}$ by:
$${}^x j(s) = j(s) \cdot \{ \varpi, x(s) \}.$$
The reduction of this section is the twist of $\alg{\bar j}$ by $x$, as desired.
\qed

\begin{corollary}
\label{SplitWish}
Let $\{y_1, \ldots, y_n \}$ be a basis of the free $\ZZ$-module $Y$, and $\alg{\bar y}_i: \alg{\bar G}_{\mult} \rightarrow \alg{\bar S}$ the cocharacter over $\FF$ corresponding to $y_i$.  Let $\alg{\bar y}_i': \alg{\bar G}_{\mult} \rightarrow \alg{\bar S}'$ be any cocharacters of $\alg{\bar S}'$ lifting the cocharacters $\alg{\bar y}_i$.  Then, there exists a section $\alg{j}$ of the pointed $\alg{K}_2$-torsor, whose reduction $\alg{\bar j}$ satisfies
$$\alg{\bar y}_i' = \alg{\bar j} \circ \alg{\bar y}_i, \mbox{ for all } 1 \leq i \leq n.$$
\end{corollary}
\proof
For every $1 \leq i \leq n$, there exists an integer $\lambda_i$ such that
$$\alg{\bar y}_i' = [\alg{\bar j} \circ \alg{\bar y}_i] \cdot \alg{\lambda}_i,$$
where $\alg{\lambda}_i$ denotes the element of $\Hom(\alg{G}_{\mult}, \alg{G}_{\mult})$ corresponding to the integer $\lambda_i$.  The collection of integers $\lambda_i$ may be assembled into an element $\lambda$ of $X = \Hom(Y, \ZZ)$ satisfying $\lambda(y_i) = \lambda_i$, since $\{ y_1, \ldots, y_n \}$ is a basis of the free $\ZZ$-module $Y$.  By twisting $\alg{\bar j}$ by $-\lambda$, as in the previous proposition, the corollary is proven.
\qed

\subsection{Parahorics}
Assume now that the residue field $\FF = \OO / \varpi \OO$ is algebraically closed (and $F$ is discretely valued as before).  The case of a (quasi-)finite residue field will follow later from \'etale descent.  Let $\alg{G}$ be a connected reductive group over $F$; thus $\alg{G}$ is quasisplit over $F$.  Let $G = \alg{G}(F)$, and let $\Build$ be the (enlarged) building of $G$ over $F$.  For $x \in \Build$, we write $G_x$ for the parahoric subgroup, which is contained in the isotropy group of $G$ fixing $x$.   

From Bruhat-Tits (Theorem 3.8.3, see also Section 4.6.2 of \cite{BT2}), there is a smooth group scheme $\alg{\underline{G}}_x$ over $\OO$, uniquely determined up to unique isomorphism with the following properties:
\begin{itemize}
\item
The generic fibre of $\underline{\alg{G}}_x$ is equal to $\alg{G}$ as group schemes over $F$.
\item
The $\OO$-points $\underline{\alg{G}}_x(\OO)$ are equal to $G_x$ as a subset of $G = \alg{\underline{G}}_x(F) = \alg{G}(F)$.  
\end{itemize}
\item
We follow the ``connected special fibre'' convention for parahoric subgroups:  the special fibre $\alg{\bar G}_x$ is a connected group scheme over $\FF$.  Let $\alg{G}'$ be a central extension of $\alg{G}$ by $\alg{K}_2$.  Let $\tilde G$ be the resulting tame extension of $G$ by $\FF^\times$.  From Construction 12.11 of \cite{B-D}, described earlier, there exists a central extension $({\alg{\bar G}_x'}, \alg{\bar p}, \alg{\bar \iota})$ of $\alg{\bar G}_x$ by $\alg{\bar G}_{\mult}$ over $\FF$, and a commutative diagram with exact rows:
$$\xymatrix{
1 \ar[r] &  \FF^\times \ar[r] \ar[d]^{=} & \tilde G_x \ar[r] \ar[d] & G_x \ar[d] \ar[r] & 1 \\
1 \ar[r] & \alg{\bar G}_{\mult}(\FF) \ar[r]^{\bar \iota} & \alg{\bar G}_x'(\FF) \ar[r]^{\bar p} & \alg{\bar G}_x(\FF) \ar[r] & 1
}$$
such that the $\tilde G_x$ is obtained via pullback from the central extension in the bottom row. 

Let $\alg{S}$ be a maximal $F$-split torus in $\alg{G}$, with canonical model $\alg{\underline{S}}$ over $\OO$.  Let $\alg{S}'$ be the resulting central extension of $\alg{S}$ by $\alg{K}_2$.  Suppose that $x$ is contained in the apartment $\Apart(S)$ of $\Build$ associated to $S$.  Then the special fibre $\alg{\bar S}$ is a maximal torus in the group $\alg{\bar G}_x$.  Letting $\alg{\bar S}' = \alg{p}^{-1}(\alg{\bar S})$, we have an extension of tori:
$$1 \rightarrow \alg{\bar G}_{\mult} \rightarrow \alg{\bar S}' \rightarrow \alg{\bar S} \rightarrow 1.$$
The above sequence of tori corresponds to a sequence of $\ZZ$-modules:
$$0 \rightarrow \ZZ \rightarrow Y' \rightarrow Y \rightarrow 0,$$
where $Y$ coincides with the cocharacter lattice of $\alg{S}$ or of $\alg{\bar S}$.

The connected group $\alg{\bar G}_x$ has a unique Levi subgroup $\alg{\bar M}_x$ containing $\alg{\bar S}$, and the connected group $\alg{\bar G}_x'$ has a unique Levi subgroup $\alg{\bar M}_x'$ containing $\alg{\bar S}'$.  In this way, we have a central extension of a reductive group by $\alg{\bar G}_{\mult}$, defined over $\FF$:
$$1 \rightarrow \alg{\bar G}_{\mult} \rightarrow \alg{\bar M}_x' \rightarrow \alg{\bar M}_x \rightarrow 1.$$

To understand this central extension of $\alg{\bar M}_x$ by $\alg{G}_{\mult}$, we apply Theorem \ref{CEByGm}, and study the resulting diagram of cocharacter lattices:
$$\xymatrix{
0 \ar[r] & \ZZ \ar[r] \ar[d]^= & Y_{x,\sconn} \times \ZZ \ar[r] \ar[d]^{\phi} & Y_{x,\sconn} \ar[r]  \ar[d]^{f_\ast} & 0 \\
0 \ar[r] & \ZZ \ar[r] & Y' \ar[r] & Y \ar[r] & 0.
}$$
Here $\alg{\bar f}: (\alg{\bar M}_x)_{\sconn} \rightarrow \alg{\bar M}_x$ is the homomorphism from the universal cover of the derived subgroup of $\alg{\bar M}_x$.  

Choosing a section $\alg{j}$ of the $\alg{K}_2$-torsor $\alg{S}' \rightarrow \alg{S}$, we get a splitting $\alg{\bar j}$ of the multiplicative $\alg{G}_{\mult}$-torsor $\alg{\bar S}' \rightarrow \alg{\bar S}$.  This, in turn, splits the sequence of cocharacter groups:
$$\xymatrix{ 0 \ar[r] & \ZZ \ar[r] & Y' \ar[r] & Y \ar[r] \ar@/_/[l]_{\bar j_\ast} & 0.}$$ 
Choosing a different section $\alg{j}$, as we saw before, will twist the splitting by an element of $X = \Hom(Y, \ZZ)$.

Such a section $\alg{j}$ yields an isomorphism $Y' \isom Y \times \ZZ$; thus it is left to determine the homomorphism $\phi_x$ fitting into the commutative diagram with exact rows:
$$\xymatrix{
0 \ar[r] & \ZZ \ar[r] \ar[d]^= & Y_{x,\sconn} \times \ZZ \ar[r] \ar[d]^{\phi_x} & Y_{x,\sconn} \ar[r]  \ar[d]^{\bar f_\ast} & 0 \\
0 \ar[r] & \ZZ \ar[r] & Y \times \ZZ \ar[r] & Y \ar[r] & 0.
}$$
\begin{remark}
The $\ZZ$-modules $Y$ and $Y'$ depend only on the torus $\alg{S}$, and not on the point $x$ in the apartment associated to $\alg{S}$.  Similarly $\alg{j}$ is chosen without regard to the point $x$, and so the splitting $\alg{\bar j}$ is also independent of the point $x$.  However, the reductive group $\alg{\bar M}_x$ and $\alg{\bar M}_{x,\sconn}$ depend on the point $x$; thus the $\ZZ$-module $Y_{x, \sconn}$ depends on the point $x$ in the apartment, hence the subscript.
\end{remark}

Let $(X, \Phi_x, Y, \Phi_x^\vee)$ denote the root datum associated to the split group $\alg{\bar M}_x$ and torus $\alg{\bar S}$ over the algebraically closed field $\FF$.  Thus $Y_{x, \sconn}$ is the $\ZZ$-submodule  of $Y$ spanned by the coroots $\Phi_x^\vee$, and the map $\bar f_\ast$ can be viewed as the inclusion.  Recall from Section 3.5.1 of \cite{Tit} the following description of this root datum:
\begin{proposition}
Let $\Phi_{\Aff, x}$ denote the set of affine roots which vanish at the point $x$.  Then $\Phi_x$ is the set of ``vector parts'' of $\Phi_{\Aff, x}$ -- the set of vector parts of affine roots vanishing at $x$.  If $\alpha \in \Phi_x$, then the coroot associated to $\alpha$ in $\Phi_x^\vee$ is equal to the coroot associated to $\alpha$ in $\Phi^\vee$, unambiguously called $\alpha^\vee$.
\end{proposition}

Let $(X', \Phi_x', Y', (\Phi_x')^\vee)$ denote the root datum associated to the group $\alg{\bar M}_x'$ and the torus $\alg{\bar S}' $.  Let $\zeta$ denote the image of $1 \in \ZZ$, under the inclusion $\ZZ \rightarrow Y'$.  The splitting $Y' = Y \times \ZZ$ (depending on $\alg{j}$, chosen independently of $x$) means that every element of $Y'$ can be written as $y + k \zeta$ for some integer $k$.  Explicitly, $\zeta$ is the inclusion $\alg{\bar \iota}$ of $\alg{\bar G}_{\mult}$ into the center of $\alg{\bar M}_x'$. 

The roots for $\alg{\bar M}_x'$ are the pullbacks of the roots of $\alg{\bar M}_x$, under the canonical homomorphism $X \rightarrow X'$.  If $\alpha \in \Phi_x$, we write $\alpha'$ for its image in $\Phi_x'$.  The coroots carry more significant information; namely, for each coroot $\alpha^\vee \in \Phi_x^\vee$, there exists a unique integer $\kappa_x(\alpha^\vee)$ satisfying 
$$(\alpha')^\vee = \alpha^\vee + \kappa_x(\alpha^\vee) \cdot \zeta.$$
The integers $\kappa_x(\alpha^\vee)$ determine the homomorphism $\phi_x$, since $\phi_x(\alpha^\vee, 0) = \alpha^\vee + \kappa_x(\alpha^\vee) \cdot \zeta$ and $Y_{x,\sconn}$ is generated by the coroots $\alpha^\vee \in \Phi_x^\vee$.

Thus by Theorem \ref{CEByGm}, these integers $\kappa_x(\alpha^\vee)$ determine the central extension $\alg{\bar M}_x'$ of $\alg{\bar M}_x$ by $\alg{\bar G}_{\mult}$, up to unique isomorphism.  Finally, we observe that if $\Delta_x = \{ \alpha_1, \ldots, \alpha_n \}$ is a system of simple roots in $\Phi_x$, then the function $\kappa_x$ is uniquely determined by its values on $\Delta_x$.  Indeed, to know a root datum, it suffices to know the character and cocharacter lattices, the roots, and the coroots associated to a system of simple roots.  The other coroots can be obtained by Weyl group reflections.

Below we summarize our approach, step-by-step:
\begin{enumerate}
\item
We wish to understand the central extension of a parahoric subgroup:
$$1 \rightarrow \FF^\times \rightarrow \tilde G_x \rightarrow G_x \rightarrow 1.$$
\item
This extension arises as the pullback of the points of a central extension of groups over the residue field:
$$1 \rightarrow \alg{\bar G}_{\mult} \rightarrow \alg{\bar G}_x' \rightarrow \alg{\bar G}_x \rightarrow 1.$$
\item
This extension arises canonically from a central extension of a Levi factor of $\alg{\bar G}_x$:
$$1 \rightarrow \alg{\bar G}_{\mult} \rightarrow \alg{\bar M}_x' \rightarrow \alg{\bar M}_x \rightarrow 1.$$
\item
Such a central extension is classified up to unique isomorphism by Theorem \ref{CEByGm}, and after choice of section $\alg{j}$ of $\alg{S}' \rightarrow \alg{S}$, is determined a single homomorphism
$$\phi_x: Y_{x, \sconn} \rightarrow Y \times \ZZ.$$
\item
To determine this homomorphism, it suffices to determine the integers $\kappa_x(\alpha^\vee)$ satisfying
$$\phi_x(\alpha^\vee) = \alpha^\vee + \kappa_x(\alpha^\vee) \cdot \zeta,$$ 
for the coroots $\alpha^\vee \in Y_{x, \sconn}$.  It even suffices to know $\kappa_x$ for the coroots of a system of simple roots in $\Phi_x$.  
\end{enumerate}

\subsection{The case $\alg{SL}_2$}
Suppose that $\alpha \in \Phi_x$, $\alpha$ is indivisible in $\Phi$, and $2 \alpha \not \in \Phi$.  Then there is a finite separable extension $E / F$, and a homomorphism with finite kernel 
$$\alg{\phi}_\alpha: \alg{R}_{E/F} \alg{SL}_{2, E} \rightarrow \alg{G},$$
with $\alg{e}_{\pm \alpha}: \alg{R}_{E/F} \alg{G}_{a,E} \rightarrow \alg{U}_{\pm \alpha} \subset \alg{G}$ as before.  These factor through the simply connected group $\alg{G}_{\sconn}$.

Define $\alg{e}_\alpha': \alg{R}_{E/F} \alg{G}_{a,E} \rightarrow \alg{G}'$ to be the canonical lift of $\alg{e}_\alpha$.  Define $\alg{n}_\alpha'$ and $\alg{h}_\alpha'$ via $\alg{e}_\alpha'$, using the same formulae used to define $\tilde n_\alpha$ and $\tilde h_\alpha$ via $\tilde e_\alpha$, so that $\alg{h}_\alpha'$ is an algebraic map (of Zariski sheaves) from $\alg{R}_{E/F} \alg{G}_{\mult}$ to $\alg{S}'$, which lifts the homomorphism $\alg{h}_\alpha: \alg{R}_{E/F} \alg{G}_{\mult} \rightarrow \alg{S}$.  This restricts to a map from $\alg{G}_{\mult}$ to $\alg{S}'$, via the natural embedding $\alg{G}_{\mult} \subset \alg{R}_{E/F} \alg{G}_{\mult}$.  This homomorphism reduces, using the residue map in K-theory, to a {\em homomorphism} which lifts the coroot $\alpha^\vee$: 
$$\alg{\bar h}_\alpha': \alg{\bar G}_{\mult}  \rightarrow \alg{\bar S}'.$$
We find that, for all $z \in \FF^\times$,
$$\bar h_\alpha'(z) = \alpha^\vee(z) \cdot \zeta(z)^{\lambda(\alpha^\vee)},$$
for some integer $\lambda(\alpha^\vee)$ depending on $\alpha$ as well as the splitting $\alg{j}$ chosen earlier.  In other words,
$$\bar h_\alpha' = \alpha^\vee + \lambda(\alpha^\vee) \cdot \zeta \in Y \times \ZZ.$$

Let $\pm a$ denote the affine roots vanishing at $x$, with vector parts $\pm \alpha$.  The affine roots $\pm a$ determine (see 1.4 of \cite{Tit} and Bruhat-Tits \cite{BT2}) subgroups $U_{\pm a}$ of the root subgroups $\alg{U}_{\pm \alpha}$, such that $U_{\pm a} \subset G_x$.  These determine an integer $m = m(a,x)$ such that:
$$U_{\add} = e_\alpha( \varpi_E^m \OO_E), \quad U_{-a} = e_{-\alpha}( \varpi_E^{-m} \OO_E),$$
where $\varpi_E$ is a uniformizing element of $E$, and $\OO_E$ the valuation ring of $E$.  

Define an element of the parahoric subgroup $G_x$ by:
$$n(a,x) = e_\alpha(\varpi_E^m) e_{-\alpha}(- \varpi_E^{-m}) e_\alpha(\varpi_E^m).$$
Then the reduction of $n(a,x)$ in $\bar G_x$ represents the Weyl reflection associated to the coroot $\alpha^\vee \in \Phi_x^\vee$.

Similarly, if we define an element of $\tilde G_x$ by:
$$\tilde n(a,x) = \tilde e_\alpha(\varpi_E^m) \tilde e_{-\alpha}(- \varpi_E^{-m}) \tilde e_\alpha(\varpi_E^m),$$
using our unipotent splitting, then $\tilde n(a,x) \in \tilde G_x$ projects onto a representative $\bar n'(a,x) \in \bar G_x'$ for the Weyl reflection associated to the root $\alpha'$.  

The formulas of Corollary \ref{WeylCover} and Example \ref{RSExamp} imply that, for all $z \in \OO_F^\times$,
\begin{eqnarray*}
\Int(\tilde n(a,x)) \tilde h_\alpha(z) & = & \tilde h_\alpha(z^{-1}) \cdot \sigma_\alpha(z^{-1}, \varpi_E^{2m}), \\
& = & \tilde h(z)^{-1} \cdot \sigma(\varpi_E^m, z)^2 \\
& = & \tilde h(z)^{-1} \cdot \{ N_{E/F} \varpi_E, z \}_{\tame}^{2qm} \\
& = & \tilde h(z)^{-1} \cdot \{ \varpi, z \}_{\tame}^{2 q m}.
\end{eqnarray*}
Here the integer $q = Q(\alpha^\vee) / [E : F]$, where $Q$ is the quadratic form associated to the central extension $\alg{G}'$ of $\alg{G}$ by $\alg{K}_2$, and $\alpha^\vee$ is viewed as a cocharacter of a maximal torus of $\alg{G}$ containing $\alg{S}$.  We also use the fact that $E/F$ is totally ramified (since $\FF$ is assumed algebraically closed), so $N_{E/F} \varpi_E$ is a uniformizing element $\varpi$ of $F$.

Reducing implies that
$$\Int(\bar n) [\alpha^\vee + \lambda \zeta](\bar z) = [-\alpha^\vee + (2mq - \lambda) \zeta](\bar z),$$
where $\lambda = \lambda(\alpha^\vee) \in \ZZ$.  On the other hand,
\begin{eqnarray*}
\Int(\bar n) [\alpha^\vee + \lambda \zeta](\bar z) & = & [s_{(\alpha')^\vee}(\alpha^\vee + \lambda \zeta)] (\bar z), \\
& = & [\alpha^\vee + \lambda \zeta - \langle \alpha^\vee + \lambda \zeta, \alpha \rangle \cdot (\alpha^\vee + \kappa_x \zeta)](\bar z), \\
& = & [\alpha^\vee + \lambda \zeta - 2 \alpha^\vee - 2 \kappa_x \zeta](\bar z), \\
& = & [- \alpha^\vee + (\lambda - 2 \kappa_x) \zeta](\bar z).
\end{eqnarray*}
It follows that 
$$(2mq - \lambda) = (\lambda - 2 \kappa_x).$$
We arrive at a fundamental relationship between integers:
\begin{equation}
\label{FormSL2}
\kappa_x(\alpha^\vee) = \lambda(\alpha^\vee) - q \cdot m(a,x).
\end{equation}
The constant $\lambda(\alpha^\vee)$ depends only on the splitting $\alg{j}$, and not on the point $x \in \Apart(S)$.  Hence the integers $\kappa_x(\alpha^\vee)$ can be computed, as $x$ varies within the apartment, from the integers $m(a,x)$ determined by the valuations on root subgroups at $x$.

\subsection{The case $\alg{SU}_3$}
Suppose that $\alpha \in \Phi_x$, $\alpha$ is indivisible in $\Phi$, and $2 \alpha \in \Phi$.  Then there is a finite separable extension $L / F$, a quadratic separable extension $E/L$, and a homomorphism with finite kernel 
$$\alg{\phi}_\alpha: \alg{R}_{L/F} \alg{SU}_{3, E/L} \rightarrow \alg{G},$$
with $\alg{e}_{\pm \alpha}: \alg{R}_{E/F} \alg{G}_{a,E} \rightarrow \alg{U}_{\pm \alpha} \subset \alg{G}$ as before.  Let $\varpi_E$ denote a uniformizing element of $E$, and let $\varpi_L = N_{E/L} \varpi_E$ and $\varpi = N_{L/F} \varpi_L$.  Write $\alg{e}_{\pm 2 \alpha}: \alg{R}_{L/F} \alg{G}_{a,L} \rightarrow \alg{G}$ also as before.  This requires us to choose a nonzero $\theta \in E$ such that $\theta + \theta^\sigma = 0$.  These homomorphisms $\alg{e}_{\pm \alpha}$, $\alg{e}_{\pm 2 \alpha}$ factor through the simply connected group $\alg{G}_{\sconn}$.

Define $\alg{e}_\alpha': \alg{R}_{L/F} \alg{J}_{E/L} \rightarrow \alg{G}'$ to be the canonical lift of $\alg{e}_\alpha$.  Define $\alg{e}_{2 \alpha}'$ to be the canonical lift of $\alg{e}_{2 \alpha}$ in the same way.  Define $\alg{n}_\alpha'$, $\alg{n}_{2 \alpha}'$, and $\alg{h}_{2 \alpha}'$ via $\alg{e}_\alpha'$ and $\alg{e}_{2 \alpha}'$, using the same formulae used to define $\tilde n_\alpha$, $\tilde n_{2 \alpha}$, and $\tilde h_{2 \alpha}$ via $\tilde e_\alpha$ and $\tilde e_{2 \alpha}$. 

Thus $\alg{h}_{2 \alpha}'$ is an algebraic map (of Zariski sheaves) from $\alg{R}_{L/F} \alg{G}_{m,L}$ to $\alg{S}'$, which lifts the homomorphism $\alg{h}_{2 \alpha}: \alg{R}_{L/F} \alg{G}_{\mult} \rightarrow \alg{S}$.  This restricts to a map from $\alg{G}_{\mult}$ to $\alg{S}'$, via the natural embedding $\alg{G}_{\mult} \subset \alg{R}_{L/F} \alg{G}_{m,L}$.  This homomorphism reduces, using the residue map in K-theory, to a {\em homomorphism} which lifts the coroot $\alpha^\vee$: 
$$\alg{\bar h}_\alpha': \alg{\bar G}_{\mult}  \rightarrow \alg{\bar S}'.$$
We find that, for all $z \in \FF^\times$,
$$\bar h_\alpha'(z) = \alpha^\vee(z) \cdot \zeta(z)^{\lambda(\alpha^\vee)},$$
for some integer $\lambda(\alpha^\vee)$ depending on $\alpha$ as well as the splitting $\alg{j}$ chosen earlier.  In other words,
$$\bar h_\alpha' = \alpha^\vee + \lambda(\alpha^\vee) \cdot \zeta \in Y \oplus \ZZ = Y + \ZZ \zeta.$$

Let $\pm a$ denote the affine roots vanishing at $x$, with vector parts $\pm \alpha$ or $\pm 2 \alpha$.  We recall from Example 1.15 of \cite{Tit} that in this situation (related to a ramified special unitary group in three variables), such vertices $x$ belong to two (mutually exclusive) types:
\begin{description}
\item[Type 1]  We say that $x$ has Type 1, if there exists an affine root $a$ vanishing at $x$, with vector part $\alpha$.  In this case, $\pm a$ determines a filtration of the root subgroup $\alg{U}_{\pm \alpha}$ of the form
$$U_{\add} = \{ e_\alpha(c,d) : \val(d) \geq \mu \},$$
where $\mu \in \val(L^\times) + \delta$.  Here, we note that $d \in E$ and $\val(E^\times) = \frac{1}{2} \val(L^\times)$, and $\delta$ is the constant:
$$\delta = \supremum \{ \val(d) : d \in E \mbox{ and } d + d^\sigma + 1 = 0 \}.$$
In odd residue characteristic, $\delta = 0$, and in even residue characteristic, $\delta < 0$.
\item[Type 2]  We say that $x$ has Type 2, if there exists an affine root $a$ vanishing at $x$, with vector part $2 \alpha$.  In this case, $\pm a$ determines a filtration of the root subgroup $\alg{U}_{\pm 2 \alpha}$ of the fom
$$U_{\add} = \{ e_{2 \alpha}(\ell) : \val(\ell \theta) \geq \mu \}.$$
\end{description}
Note that in both cases, we find a rational number $\mu = \mu(a,x) \in \val(E^\times)$.  There is an integer $m = m(a,x)$ such that $\mu = m \cdot \val(\varpi_E)$. 

Let $x$ be a vertex in $\Apart$ at which $\pm a$ vanishes.  The reduction of the parahoric $\alg{\bar G}_x$ has $\pm \alpha$ as roots in Type 1 and $\pm 2 \alpha$ as roots in Type 2, with respect to the maximal torus $\alg{\bar S}$.  We may define elements of the parahoric $G_x$ by:
$$n(a,x) = n_\alpha(c,d) \mbox{ or } n(a,x) = n_{2 \alpha}(d \theta^{-1}) = n_\alpha(0, d),$$
in Type 1 or Type 2 respectively, requiring $\val(d) = \mu = \mu(a,x)$.  Similarly, we define $\tilde n(a,x) = \tilde n_\alpha(c,d)$ or $\tilde n(a,x) = \tilde n_\alpha(0, d)$; these are elements of $\tilde G_x$, projecting onto a representative $\bar n'(a,x) \in \bar G_x'$ for the Weyl reflection associated to the root $\alpha'$.  

Define $t = \val(\theta) \cdot [E:F]$.  The formulas of Corollary \ref{WeylCover} and Example \ref{RSExamp} imply that, for all $z \in \OO_F^\times$,
\begin{eqnarray*}
\Int(\tilde n(a,x)) \tilde h_{2 \alpha}(z) & = & \tilde h_{2 \alpha}(z^{-1}) \cdot \sigma_{2 \alpha}(z^{-1}, \frac{d d^\sigma}{\theta \theta^\sigma}), \\
& = & \tilde h(z)^{-1} \cdot \sigma_{2 \alpha}(N_{E/L}(d \theta^{-1}), z) \\
& = & \tilde h(z)^{-1} \cdot \{ N_{E/F}(d \theta^{-1}), z \}_{\tame}^{q} \\
& = & \tilde h(z)^{-1} \cdot \{ \varpi, z \}_{\tame}^{2 q(m-t)}.
\end{eqnarray*}
Here the integer $q = Q(\alpha^\vee) / [L : F]$, where $Q$ is the quadratic form associated to the central extension $\alg{G}'$ of $\alg{G}$ by $\alg{K}_2$, and $\alpha^\vee$ is viewed as a cocharacter of a maximal torus of $\alg{G}$ containing $\alg{S}$.

Reducing implies that
$$\Int(\bar n) [\alpha^\vee + \lambda \zeta](\bar z) = [-\alpha^\vee + (2(m-t)q - \lambda) \zeta](\bar z),$$
where $\lambda = \lambda(\alpha^\vee) \in \ZZ$.  On the other hand, as in the $\alg{SL}_2$ case,
\begin{eqnarray*}
\Int(\bar n) [\alpha^\vee + \lambda \zeta](\bar z) & = & [s_{(\alpha')^\vee}(\alpha^\vee + \lambda \zeta)] (\bar z), \\
& = & [\alpha^\vee + \lambda \zeta - \langle \alpha^\vee + \lambda \zeta, \alpha \rangle \cdot (\alpha^\vee + \kappa_x \zeta)](\bar z), \\
& = & [\alpha^\vee + \lambda \zeta - 2 \alpha^\vee - 2 \kappa_x \zeta](\bar z), \\
& = & [- \alpha^\vee + (\lambda - 2 \kappa_x) \zeta](\bar z).
\end{eqnarray*}
It follows that 
$$(2(m-t)q - \lambda) = (\lambda - 2 \kappa_x).$$
We arrive at a fundamental relationship:
\begin{equation}
\label{FormSU3}
\kappa_x(\alpha^\vee) = \lambda(\alpha^\vee) - q \cdot (m(a,x) - t).
\end{equation}
The constant $\lambda(\alpha^\vee)$ depends only on the splitting $\alg{j}$, and not on the point $x \in \Apart(S)$.  The constant $t = \val(\theta) \cdot [E:F]$ can often be chosen to be $0$ or $1$, in in any case does not depend on the point $x$.  Hence the integers $\kappa_x(\alpha^\vee)$ can be computed, as $x$ varies within the apartment, from the integers $m(a,x)$ determined by the valuations on root subgroups at $x$.

\subsection{The simplest example}

Consider the simplest example, $\alg{G} = \alg{SL}_2$, and $\alg{G}'$ the central extension of $\alg{G}$ by $\alg{K}_2$ associated to the quadratic form $Q = Q_1$.  Let $\pm \alpha$ denote the roots of $\alg{G}$ with respect to the usual torus $\alg{S}$ of diagonal matrices.  Thus $Q(\alpha^\vee) = 1$.  Fix the usual Chevalley system on $\alg{G}$, yielding a base point $x_0$ in the apartment $\Apart = \Apart(S)$ in the building $\Build = \Build(G)$.  

Let $x$ be a vertex in $\Apart$; thus $x = x_0 - \frac{m}{2} \alpha^\vee$ for some integer $m$.  The parahoric subgroup at $x$ looks like:
$$G_x = \Matrix{\OO}{\ideal{p}^m}{\ideal{p}^{-m}}{\OO} \cap SL_2(F).$$
The reduction $\alg{\bar G}_x$ is isomorphic to $\alg{\overline{SL}}_2$ over $\FF$, with maximal torus $\alg{\bar S}$.  A representative for the nontrivial Weyl element in $\alg{\bar G}_x$ is obtained from the reduction of
$$n = \Matrix{0}{\varpi^m}{\varpi^{-m}}{0}.$$

At each such point $x$, the construction of Brylinski and Deligne yields a central extension:
$$1 \rightarrow \alg{\bar G}_{\mult} \rightarrow \alg{\bar G}_x' \rightarrow \alg{\bar G}_x \rightarrow 1.$$
Choose a splitting $\alg{j}$ of the resulting central extension of tori:
$$1 \rightarrow \alg{\bar G}_{\mult} \rightarrow \alg{\bar S}' \rightarrow \alg{\bar S} \rightarrow 1.$$

Let $(X', \Phi_x', Y', (\Phi_x')^\vee)$ denote the root datum of $\alg{\bar G}_x'$ with respect to $\alg{\bar S}'$.  Our splitting $\alg{j}$ identifies:
$$X' = X \oplus \ZZ = \ZZ \oplus \ZZ, \quad Y' = Y \oplus \ZZ = \ZZ \oplus \ZZ.$$
This is independent of the choice of point $x \in \Apart$.  

As $x = x_0 - \frac{m}{2}$ varies within the apartment -- as $m$ varies over integers -- the roots stay constant and the coroots for $\alg{\bar G}_x'$ vary as follows:
\begin{eqnarray*}
\Phi_x' & = & \{ \alpha + 0, -\alpha + 0 \} = \{ (2,0), (-2,0) \}, \\
(\Phi_x')^\vee & = &  \{ \alpha^\vee - m, -\alpha^\vee + m \} = \{ (1,-m), (-1,+m) \}.
\end{eqnarray*}
All the groups $\alg{\bar G}_x'$ are uniquely isomorphic to $\alg{\overline{SL}}_2 \times \alg{\bar G}_{\mult}$, but the root datum changes, with respect to their common maximal torus.  Of course, this description depends on a choice of base point $x_0$, and splitting $\alg{j}$ among other things.  But any change in these choices simply shifts $m$ by an integer.  Using a central extension $\alg{G}'$ arising from the quadratic form $Q = d \cdot Q_1$ would replace $m$ by $dm$ in the above description.

\subsection{The split group $\alg{G}_2$}

Consider $\alg{G} = \alg{G}_2$, the split Chevalley group of type $\Type{G}_2$ over $F$, with split maximal torus $\alg{S}$.  Let $\alg{G}'$ be the central extension of $\alg{G}$ by $\alg{K}_2$ associated to the quadratic form $Q = Q_1$.  Let $\alpha$ and $\beta$ be simple positive roots, with $\alpha$ short and $\beta$ long, so that the positive roots are:
$$\Phi^+ = \{ \alpha, \beta, \beta + \alpha, \beta + 2 \alpha, \beta + 3 \alpha, 2 \beta + 3 \alpha \}.$$
Note that $\alpha^\vee$ is long and $\beta^\vee$ is short, so $Q(\alpha^\vee) = 3$ and $Q(\beta^\vee) = 1$.  Note also that
$$\langle \alpha^\vee, \beta \rangle = -3, \quad \langle \beta^\vee, \alpha \rangle = -1.$$  

There are three types of vertices in the building $\Build$ of $G$, with local spherical buildings of type $\Type{G}_2$ (hyperspecial vertices), $\Type{A}_2$, and $\Type{A}_1 \times \Type{A}_1$.  At each type of vertex, we find a different group $\alg{\bar G}_x$, and we describe the central extensions arising from $\alg{G}'$ here.  We fix a section $\alg{j}$ of $\alg{S}' \rightarrow \alg{S}$ in such a way that $\lambda(\alpha^\vee) = \lambda(\beta^\vee) = 0$, since these coroots form a basis for the $\ZZ$-module $Y$ (see Corollary \ref{SplitWish}).  Thus we we identify
$$X' = X \oplus \ZZ = \ZZ \alpha + \ZZ \beta + \ZZ, \quad Y' = Y \oplus \ZZ = \ZZ \alpha^\vee + \ZZ \beta^\vee + \ZZ \zeta.$$

At the hyperspecial point $x_0$ corresponding to our initial Chevalley system, we have $\alg{\bar G}_{x_0}' = \alg{\bar G}_{x_0} \times \alg{\bar G}_{\mult} = \alg{\bar G}_2 \times \alg{\bar G}_{\mult}$.  The root datum at this base point is given by 
$$\Phi_{x_0}' = \{ (\gamma,0) \}_{\gamma \in \Phi}, \quad (\Phi_{x_0}')^\vee = \{ (\gamma^\vee,0) \}_{\gamma^\vee \in \Phi^\vee}.$$

A nearby hyperspecial point lies at $x = x_0 + \alpha^\vee + \beta^\vee$, at the intersection of two affine root hyperplanes:
$$[\alpha + \beta](x - x_0) = [\beta](x - x_0) + 1 = 0.$$
Observe that:
$$\alpha(x - x_0) = 1, \quad \beta(x - x_0) = -1.$$
Letting $a = \alpha - 1$ and $b = \beta + 1$ be the resulting affine roots, we find root subgroups of the parahoric $G_x$
$$U_{\add} = e_\alpha(\ideal{p}), \quad U_b = e_\beta(\ideal{p}^{-1}).$$

It follows that in the root datum of $\alg{\bar G}_x'$ we find roots $\alpha' = (\alpha,0)$ and $\beta' = (\beta,0)$ in $X' = X \oplus \ZZ$.  The associated coroots in $Y'$ are computed with Equation \ref{FormSL2}:
$$(\alpha')^\vee = (\alpha^\vee, -3), \quad (\beta')^\vee = (\beta^\vee, 1).$$
Since the coroots are coplanar in the vector space $Y' \otimes_\ZZ \QQ \isom \QQ^3$, we can find all the coroots in $(\Phi_x')^\vee$ from these two.
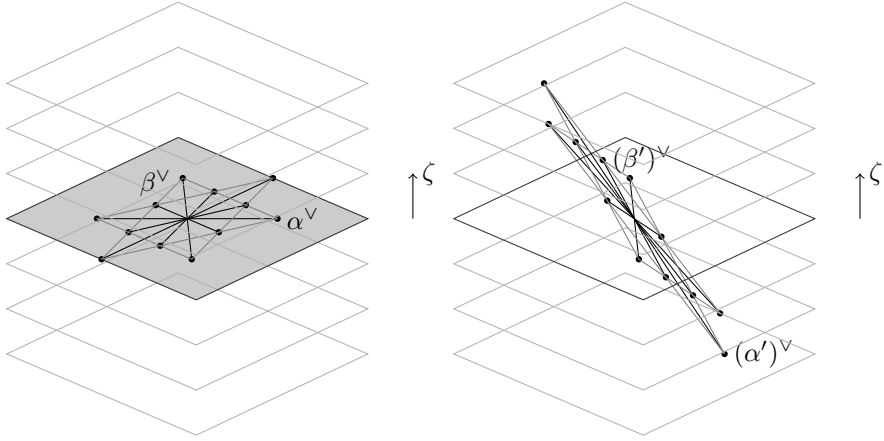
\begin{figure}
\begin{tikzpicture}[x={(2cm,0cm)},y={(-0.7cm,0.3cm)},z={(0cm,1cm)}, scale=0.6]
\filldraw [black!20] (-2,-6,0) -- (2,0,0) -- (2,6,0) -- (-2,0,0) -- (-2,-6,0);
\foreach \x/\y in {1/0, 0/1, 1/1, 1/2, 1/3, 2/3}
  { \draw (0,0) -- (\x, \y, 0); 
   \draw (0,0) -- (-\x, -\y, 0);
   \fill (\x, \y) circle (2pt);
   \fill (-\x, -\y) circle (2pt);
   }
   \draw[black!50] (1,0) -- (1,3) -- (-2,-3) -- (1,0);
   \draw[black!50] (2,3) -- (-1,0) -- (-1,-3) -- (2,3);
   \draw (1,0) node[right] {$\alpha^\vee$};
   \draw (0,1) node[above] {$\beta^\vee$};
   
   \foreach \z/\c in {-3/30,-2/30,-1/30, 0/80, 1/30, 2/30,3/30}
      \draw[black!\c] (-2,-6,\z) -- (2,0,\z) -- (2,6,\z) -- (-2,0,\z) -- (-2,-6,\z);

   \draw[->] (2.5,0,0) -- (2.5,0,1) node[right] {$\zeta$};
\end{tikzpicture}
\begin{tikzpicture}[x={(2cm,-3cm)},y={(-0.7cm,1.3cm)},z={(0cm,1cm)}, scale=0.6]
\foreach \x/\y in {1/0, 0/1, 1/1, 1/2, 1/3, 2/3}
  { \draw (0,0) -- (\x, \y, 0); 
   \draw (0,0) -- (-\x, -\y, 0);
   \fill (\x, \y) circle (2pt);
   \fill (-\x, -\y) circle (2pt);
   }
   \draw[black!50] (1,0) -- (1,3) -- (-2,-3) -- (1,0);
   \draw[black!50] (2,3) -- (-1,0) -- (-1,-3) -- (2,3);
   \draw (1,0) node[right] {$(\alpha')^\vee$};
   \draw (0,1) node[right] {$(\beta')^\vee$};
   \foreach \z/\c in {-3/30,-2/30,-1/30, 0/80, 1/30, 2/30,3/30}
      \draw[black!\c, x={(2cm,0cm)},y={(-0.7cm,0.3cm)},z={(0cm,1cm)}] (-2,-6,\z) -- (2,0,\z) -- (2,6,\z) -- (-2,0,\z) -- (-2,-6,\z);

   \draw[->] (2.5,0,7.5) -- (2.5,0,8.5) node[right] {$\zeta$};
\end{tikzpicture}
\caption{The coroots $(\Phi_{x_0}')^\vee$ on the left, and the coroots $(\Phi_x')^\vee$ on the right.  Both are embedded in the same $\ZZ$-module $Y' = Y + \zeta \ZZ$.  The shaded plane containing the coroots has been skewed, sending $\alpha^\vee$ to $\alpha^\vee - 3 \zeta$ and $\beta^\vee$ to $\beta^\vee + \zeta$.} 
\end{figure}

Now consider a point $y \in \Apart$ at which the local Dynkin diagram has type $\Type{A}_1 \times \Type{A}_1$.  Such a point occurs at the midpoint of the segment from $x_0$ to $x$; this midpoint is $y = x_0 + \frac{1}{2}(\alpha^\vee + \beta^\vee)$.  Then $\alg{\bar G}_y$ is a group isomorphic to $\alg{\overline{SO}}_4$, which is neither simply-connected nor adjoint.

The only roots at $y$ (the vector parts of affine roots vanishing at $y$) are the following:
$$\Phi_y = \{ \pm (\alpha + \beta), \pm (3 \alpha + \beta) \}.$$
Let $\gamma = \alpha + \beta$ and $\delta = 3 \alpha + \beta$, so $\Phi_y = \{ \pm \gamma, \pm \delta \}$.  The associated coroots are:
$$\gamma^\vee = \alpha^\vee + 3 \beta^\vee, \quad \delta^\vee = \alpha^\vee + \beta^\vee.$$

Then we find that
$$\gamma(y - x_0) = 0, \quad \delta(y - x_0) = 1.$$
Let $c = \gamma + 0$ and $d = \delta = d - 1$ be the associated affine roots vanishing at $y$.  There are corresponding subgroups of the parahoric $G_y$:
$$U_c = e_\gamma(\OO), \quad U_d = e_\delta(\ideal{p}).$$
The associated coroots are
$$(\gamma')^\vee = (\gamma^\vee, 0), \quad (\delta')^\vee = (\delta^\vee, -1),$$
using the fact that $\delta^\vee$ is a short coroot and so $Q(\delta^\vee) = 1$.
\begin{figure}
\begin{tikzpicture}[x={(2cm,0cm)},y={(-0.7cm,0.3cm)},z={(0cm,1cm)}, scale=0.6]
\filldraw [black!20] (-2,-6,0) -- (2,0,0) -- (2,6,0) -- (-2,0,0) -- (-2,-6,0);
\foreach \x/\y in {1/3, 1/1}
  { \draw (0,0) -- (\x, \y, 0); 
   \draw (0,0) -- (-\x, -\y, 0);
   \fill (\x, \y) circle (2pt);
   \fill (-\x, -\y) circle (2pt);
   }
   \fill (1,3) circle(2pt);
   \fill (1,1) circle(2pt);
   \draw[black!50] (1,0) -- (1,3) -- (-2,-3) -- (1,0);
   \draw[black!50] (2,3) -- (-1,0) -- (-1,-3) -- (2,3);
   \draw (1,3) node[above] {$\gamma^\vee$};
   \draw (1,1) node[above] {$\delta^\vee$};
   \foreach \z/\c in {-3/30,-2/30,-1/30, 0/80, 1/30, 2/30,3/30}
      \draw[black!\c] (-2,-6,\z) -- (2,0,\z) -- (2,6,\z) -- (-2,0,\z) -- (-2,-6,\z);

   \draw[->] (2.5,0,0) -- (2.5,0,1) node[right] {$\zeta$};
\end{tikzpicture}
\begin{tikzpicture}[x={(2cm,-1.5cm)},y={(-0.7cm,0.8cm)},z={(0cm,1cm)}, scale=0.6]
\filldraw [black!20] (-2,-6,0) -- (2,0,0) -- (2,6,0) -- (-2,0,0) -- (-2,-6,0);
\foreach \x/\y in {1/3, 1/1}
  { \draw (0,0) -- (\x, \y, 0); 
   \draw (0,0) -- (-\x, -\y, 0);
    \fill (\x, \y) circle (2pt);
    \fill (-\x, -\y) circle (2pt);
   }
   \draw[black!50] (1,0) -- (1,3) -- (-2,-3) -- (1,0);
   \draw[black!50] (2,3) -- (-1,0) -- (-1,-3) -- (2,3);
   \draw (1,3) node[above] {$(\gamma')^\vee$};
   \draw (1,1) node[above] {$(\delta')^\vee$};
   \foreach \z/\c in {-3/30,-2/30,-1/30, 0/80, 1/30, 2/30,3/30}
      \draw[black!\c, x={(2cm,0cm)},y={(-0.7cm,0.3cm)},z={(0cm,1cm)}] (-2,-6,\z) -- (2,0,\z) -- (2,6,\z) -- (-2,0,\z) -- (-2,-6,\z);

   \draw[->,  x={(2cm,0cm)},y={(-0.7cm,0.3cm)},z={(0cm,1cm)}] (2.5,0,0) -- (2.5,0,1) node[right] {$\zeta$};
\end{tikzpicture}
\caption{The coroots $(\Phi_{x_0}')^\vee$ on the left (with $\gamma^\vee + 0 \zeta$ and $\delta^\vee + 0 \zeta$ highlighted), and the coroots $(\Phi_y')^\vee$ on the right.  Both are embedded in the same $\ZZ$-module $Y' = Y + \zeta \ZZ$.  The planes containing the coroots are shaded.} 
\end{figure}
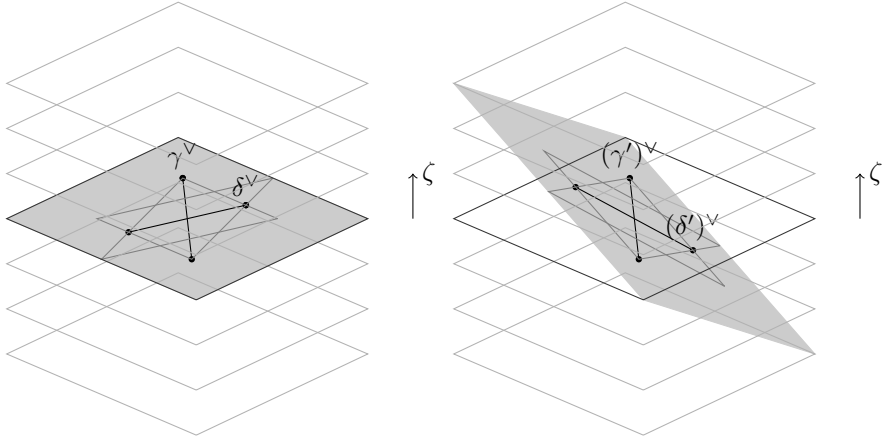

At the point $y$, with $\alg{\bar G}_y \isom \alg{\overline{SO}}_4$, we have computed the root datum of the central extension $\alg{\bar G}_y'$ -- the central extension splits and $\alg{\bar G}_y'$ is isomorphic to $\alg{\overline{SO}}_4 \times \alg{\bar G}_{\mult}$.

\bibliographystyle{plain}
\bibliography{MetaBib}
\end{document}